\documentclass{amsart}
\usepackage{amsmath,amssymb,color}
\usepackage[all]{xy}

\newcommand{\N}{{\mathbb{N}}}

\newcommand{\R}{{\mathbb{R}}}

\newcommand{\Ch}{{\mathcal C}}
\newcommand{\Dh}{{\mathcal D}}
\newcommand{\Eh}{{\mathcal E}}
\newcommand{\Fh}{{\mathcal F}}
\newcommand{\Gh}{{\mathcal G}}

\newcommand{\Mh}{{\mathcal M}}

\newcommand{\Oh}{{\mathcal O}}

\newcommand{\Zh}{{\mathcal Z}}

\newcommand{\ad}{\mathrm{ad}\,}

\newcommand{\Aut}{\mathrm{Aut}\,}
\newcommand{\be}{\mathbf{1}}

\newcommand{\id}{\mathrm{id}}

\newcommand{\pr}{\mathrm{pr}}

\newcommand{\supp}{\mathrm{supp }\,}

\newcommand{\SSS}{\scriptscriptstyle}

\newcommand{\ot}{\otimes}
\newcommand{\ep}{{\varepsilon}}

\newcommand{\U}{{\mathcal U}}
\newcommand{\fset}{{\mathcal F}}
\newcommand{\gset}{{\mathcal G}}
\newcommand{\End}{\mathrm{End}\,}

\newcounter{number}[section]

\newenvironment{nummer}{\refstepcounter{number}
{\noindent\arabic{section}.\arabic{number}}}{}

\newcommand{\bn}{\begin{nummer} \rm}
\newcommand{\en}{\end{nummer}}

\newenvironment{thms}{\noindent {\sc Theorem:} \it}{}
\newenvironment{lms}{\noindent {\sc Lemma:} \it}{}
\newenvironment{props}{\noindent {\sc Proposition:} \it}{}
\newenvironment{dfs}{\noindent {\sc Definition:} \it}{}
\newenvironment{cors}{\noindent {\sc Corollary:} \it}{}

\newenvironment{rems}{\noindent {\sc Remark:}}{}

\newenvironment{nots}{\noindent {\sc Notation:} }{}

\newenvironment{nproof}{\noindent {\sc Proof:}}{\mbox{}\hfill
\rule[-.2ex]{.25em}{1.8ex}}

\newenvironment{nnproof}{\noindent {\sc Proof}}{\mbox{}\hfill
\rule[-.2ex]{.25em}{1.8ex}}

\begin{document}

\title{\sc Trivialization of $\Ch(X)$-algebras with strongly self-absorbing fibres}

\author{Marius Dadarlat}
\address{Department of Mathematics, Purdue University,
West Lafayette,\\
IN 47907, USA}

\email{mdd@math.purdue.edu}

\author{Wilhelm Winter}
\address{Mathematisches Institut der Universit\"at M\"unster\\
Einsteinstr.\ 62\\ D-48149 M\"unster, \indent Germany}

\email{wwinter@math.uni-muenster.de}

\date{\today}
\thanks{{\it Supported by:} The first named author was partially supported by
NSF grant \#DMS-0500693. \\
\indent The second named author was supported by the DFG (SFB 478), and by the EU-Network\\
\indent  Quantum Spaces - Noncommutative Geometry (Contract No.\ HPRN-CT-2002-00280)}

\setcounter{section}{-1}

\begin{abstract}
Suppose $A$ is a separable unital $\Ch(X)$-algebra each fibre of which is isomorphic to the same strongly self-absorbing and $K_{1}$-injective $C^{*}$-algebra $\Dh$. We show that $A$ and  $\Ch(X) \otimes \Dh$ are  isomorphic as $\Ch(X)$-algebras provided the compact Hausdorff space $X$ is finite-dimensional. This statement is known not to extend to the infinite-dimensional case. 
\end{abstract}

\maketitle

\section{Introduction}

A unital and separable $C^{*}$-algebra $\Dh\neq \mathbb{C}$ is strongly self-absorbing if there is an isomorphism $\Dh\stackrel\cong\rightarrow\Dh \otimes \Dh$
 which is approximately unitarily equivalent
to the inclusion map $\Dh\to \Dh \otimes \Dh$, $d\mapsto d \otimes \be_{\Dh}$, cf.\ \cite{TomsWinter:ssa}.
 Strongly self-absorbing $C^{*}$-algebras are known to be simple and nuclear; moreover, they are either purely infinite or stably finite.
The only known examples  are UHF algebras of infinite type (i.e., every prime number that occurs in the respective supernatural number occurs with infinite multiplicity), the Cuntz algebras $\Oh_{2}$ and $\Oh_{\infty}$, the Jiang--Su algebra $\Zh$ and tensor products of $\Oh_{\infty}$ with UHF algebras of infinite type, see \cite{TomsWinter:ssa}. All these examples are known to be $K_{1}$-injective,
 i.e., the canonical map $U(\Dh)/U_{0}(\Dh) \to K_{1}(\Dh)$ is injective.
  
  The notion of $\Ch(X)$-algebras, 
introduced by Kasparov, is a generalization of  continuous bundles (or fields) of
$C^{*}$-algebras,  cf.\
\cite{Kas:Novikov}.
The main result of our paper is the following:

\bn
\begin{thms}\label{thm:ssa-fields-are-trivial}\label{C(X)-triviality}
 Let $A$ be a separable unital  $\Ch(X)$-algebra
over a finite dimensional compact metrizable space  $X$. Suppose that
 all the fibres of $A$ are isomorphic to the same
 strongly self-absorbing $K_1$-injective $C^{*}$-algebra $\Dh$.
 Then, $A$  and $\Ch(X)\otimes \Dh$ are isomorphic
  as $\Ch(X)$-algebras.
\end{thms}

In the case where $\Dh=\Oh_{\infty}$ or $\Dh=\Oh_{2}$, the preceding result was already obtained by the first named author in \cite{Dad:bundles-fdspaces}; it is new for UHF algebras of infinite type and for the Jiang--Su algebra. While it would already be quite satisfactory to have this trivialization result for the known strongly self-absorbing examples, it is remarkable that it can be proven without any of the special properties of the \emph{concrete} algebras. In this sense, the theorem further illustrates the importance of strongly self-absorbing algebras for the Elliott program. (In fact, Theorem~\ref{thm:ssa-fields-are-trivial} may clearly be regarded as a classification result for $\Ch(X)$-algebras with strongly self-absorbing fibres, and as such it contributes to the non-simple case of the Elliott conjecture.) This point of view is one of the reasons why we think that not only the theorem, but also the methods developed for its proof are of independent interest.  In the subsequent sections we shall therefore present two rather different proofs of Theorem~\ref{thm:ssa-fields-are-trivial}. The first approach follows the strategy of \cite{Dad:bundles-fdspaces}, using the main results of \cite{HirshbergRordamWinter:absorb-ssa} and \cite{Dadarlat-Winter:KK_D-topology}. The second one  follows ideas from  Section 4 of \cite{HirshbergRordamWinter:absorb-ssa}. We outline both approaches in \ref{outline-marius-approach} and \ref{outline-wilhelms-approach}.
\en

In \cite[Example~4.7]{HirshbergRordamWinter:absorb-ssa}, Hirshberg, R{\o}rdam and the second named author constructed examples of  $\Ch(X)$-algebras with $X = \prod_{\N} S^{2}$  and fibres isomorphic to any prescribed UHF algebra of infinite type, which do not absorb this UHF algebra (and hence cannot be trivial).
In \cite{Dad:KKX-equiv}, the first named author modified this example to construct a separable, unital $\Ch(X)$-algebra over the Hilbert cube with each fibre isomorphic to $\Oh_{2}$, but which does not have trivial $K$-theory (and hence cannot be trivial either).
Therefore, the dimension condition on $X$ in Theorem~\ref{thm:ssa-fields-are-trivial} cannot be removed. However, at the present stage, it is not known whether the Theorem also fails for infinite dimensional spaces $X$ if $\Dh=\Oh_{\infty}$ or $\Dh=\Zh$.

\section{$\Ch(X)$-Algebras }

\noindent
We recall some facts and notation about $\Ch(X)$-algebras (introduced in \cite{Kas:Novikov}). For our purposes, it will be enough to restrict to
compact spaces.

\bn
\label{D-Ch(X)}
\begin{dfs}
Let $A$ be a $C^{*}$-algebra and $X$ a  compact
Hausdorff space.
$A$ is a $\Ch(X)$-algebra, if there is a unital $*$-homomorphism
$\mu \colon  \Ch(X) \to \Zh(\Mh(A))$ from
$\Ch(X)$ to the center of the multiplier algebra of $A$.
\end{dfs}

The map $\mu$ is called the \emph{structure map}. We will
not always write it explicitly.

If $A$ is as above and  $Y \subset X$ is a closed subset, then
\[
J_{Y}:= \Ch_{0}(X \setminus Y)  \cdot A
\]
is a (closed) two-sided ideal of $A$; we denote the quotient map by
$\pi_{Y}$ and set
\[
A(Y)=A_{Y}:= A/J_{Y}.
\]
$A_{Y}$ may be regarded as a $\Ch(X)$-algebra or as a $\Ch(Y)$-algebra in the obvious way.

If $a \in A$, we sometimes write $a_{Y}$ for $\pi_{Y}(a)$. If $Y$
consists of just one point $x$, we will slightly abuse notation
and write $A_{x}$ (or $A(x)$), $J_x$, $\pi_{x}$ and $a_{x}$ (or $a(x)$)
in place of
$A_{\{x\}}$, $J_{\{x\}}$, $\pi_{\{x\}}$ and $a_{\{x\}}$,
respectively. We say that $A_{x}$ is the fibre of $A$ at
$x$. If $\varphi:A\to B$ is a morphism of $\Ch(X)$-algebras, we denote
by $\varphi_Y$ the corresponding restriction map $A_Y\to B_Y$.
\en

\bn \label{upper-semicontinuity}
For any $a \in A$  we have
\[
\|a\| = \sup\{\|a_{x}\| \,:\, x \in X\}.
\]
Moreover, the function $x \mapsto \|a_{x}\|$ from $X$ to
$\R$ is upper semicontinuous.
If the map  is continuous for any $a \in A$, then $A$ is said
to be a continuous $\Ch(X)$-algebra. By \cite[Lemma~2.3]{Dad:bundles-fdspaces},
any unital $\Ch(X)$-algebra with simple nonzero fibres is automatically continuous. In this case the structure map is injective and hence $X$
is metrizable if $A$ is separable.
\en

\section{Proving the main result: The first approach}

In this section we give a proof of Theorem~\ref{thm:ssa-fields-are-trivial} which follows ideas of \cite{Dad:bundles-fdspaces} and relies on the main absorption
result Theorem~4.6 of \cite{HirshbergRordamWinter:absorb-ssa} and on \cite[Theorem~2.2]{Dadarlat-Winter:KK_D-topology}:

\bn
\label{D-asu}
\begin{thms}
Let $A$ and $\Dh$ be  unital $C^{*}$-algebras,
with $\Dh$ separable, strongly self-absorbing and $K_{1}$-injective.
 Then, any two unital $*$-homomorphisms $\sigma, \gamma:\Dh \to A \otimes \Dh$ are strongly asymptotically unitarily equivalent, i.e.
 there is a unitary-valued continuous map $u:(0,1]\to
( A\otimes \Dh)$, $t\mapsto u_t$, with $u_1=\be_{ A\otimes \Dh}$ and such that
$\lim_{t \to 0} \|u_t d u_t^*-\alpha(d)\|=0$ for all $ d \in
 D$.
\end{thms}
\en

\bn
\label{outline-marius-approach}
Let us give an outline of the strategy.
While a unital, separable,
 strongly self-absorbing  $C^{*}$-algebra $\Dh$ is not in general
 weakly semiprojective, we prove nevertheless that $\Dh$ satisfies
 a property analogous to weak semiprojectivity in the category of
 unital and $\Dh$-stable $C^{*}$-algebras, see Proposition~\ref{prop:wsp-for-ssa}.
 If $A$ is a unital, separable, $\Ch(X)$-algebra over a finite dimensional compact space with fibres isomorphic to $\Dh$, then $A$ is $\Dh$-stable by \cite{HirshbergRordamWinter:absorb-ssa}. This enables us to
 use the relative weak semiprojectivity of $\Dh$ described above and approximate
 $A$ by certain pullback $\Ch(X)$-subalgebras with with fibres isomorphic to $\Dh$, see Proposition~\ref{basic-approx}.
 Using Theorem~\ref{D-asu}, we can deform continuously families of endomorphisms to families of automorphisms of $\Dh$, and thus show that these pullback $\Ch(X)$-subalgebras
 are actually trivial, see Proposition~\ref{local-triviality-CX-semiprojective}. We then complete the proof by applying
 Elliott's intertwining argument.
\en

\bn
Let $\eta:B \to A$ and $\psi:E \to A$ be $*$-homomorphisms. The pullback of these
maps is
\[B\oplus_{\eta,\psi}E=\{(b,e)\in B \oplus E:\,
\eta(b)=\psi(e)\}.\] We are going to use pullbacks in the context of
$\Ch(X)$-algebras.
\en

We need the following three lemmas from \cite{Dad:bundles-fdspaces},
reproduced here for the convenience of the reader.

\bn
\begin{lms}\label{a(x)delta}
  Let $A$ be a
    $\Ch(X)$-algebra and let $B \subset A$ be a $\Ch(X)$-subalgebra.
    Let $a \in A$ and let $Y$ be a closed  subset of $X$.
    \begin{itemize}
\item[(i)] The map $x\mapsto \|a(x)\|$ is
    upper semi-continuous.
\item[(ii)] $\|\pi_Y(a)\|=\max\{\|\pi_x(a)\|:x\in Y\}$
        \item[(iii)] If $a(x)\in \pi_x(B)$ for all $x \in X$, then $a \in B$.
        \item[(iv)] If $\delta>0$ and $a(x)\in_\delta \pi_x(B)$
        for all $x \in X$, then $a \in_\delta B$.
\item[(v)] The restriction of $\pi_x:A \to A(x)$ to $B$ induces an isomorphism
$B(x)\cong\pi_x(B)$ for all $x \in X$.
    \end{itemize}
\end{lms}
\en

 Let $B\subset A(Y)$ and $E\subset A(Z)$ be
$\Ch(X)$-subalgebras such that $ \pi^Z_{Y\cap Z} (E)\subseteq \pi^Y_{Y\cap
Z}(B)$. By a basic property of $\Ch(X)$-algebras, see \cite[Lemma~2.4]{Dad:bundles-fdspaces},  the
pullback $B\oplus_{\pi^Z_{Y\cap Z}, \pi^Y_{Y\cap Z}}E$ is isomorphic to the
$\Ch(X)$-subalgebra $B\oplus_{Y\cap Z}E$ of $A$ defined as
\[B\oplus_{Y\cap Z}E=\{a\in A: \pi_Y(a)\in B, \pi_Z(a)\in E\}.\]

\bn
\begin{lms}\label{sumover-Y} The fibres of $B\oplus_{Y\cap Z}E$ are given by
\[\pi_x(B\oplus_{Y\cap Z}E)=\left\{%
\begin{array}{ll}
    {\pi_x(B)}, & \hbox{if $x\in X\setminus Z,$} \\
    {\pi_x(E)}, & \hbox{if $x\in Z,$} \\
\end{array}%
\right.\] and there is an exact sequence of $C^{*}$-algebras
\begin{equation}\label{use}
    \xymatrix{
 {0}\ar[r]& {\{b\in B:\pi_{Y\cap Z}(b)=0\}}\ar[r]
                & {B\oplus_{Y\cap Z}E} \ar[r]^{\,\,\quad\pi_Z}&{E}\ar[r]&0
}
\end{equation}
\end{lms}
\en
Let $X$, $Y$, $Z$ and $A$ be as above. Let $\eta:B\hookrightarrow A(Y)$ be
 a $\Ch(Y)$-linear $*$-monomorphism
     and let $\psi:E \hookrightarrow A(Z)$
     be a $\Ch(Z)$-linear $*$-monomorphism. Assume that
     \begin{equation}\label{inclusion}
   \pi^Z_{Y\cap Z} (\psi(E))\subseteq \pi^Y_{Y\cap Z} (\eta(B)).
\end{equation} This gives a map $\gamma=\eta_{Y\cap Z}^{-1}\psi_{Y\cap Z}: E(Y\cap Z)
\to B(Y\cap Z)$. To simplify notation we let $\pi$ stand for both $\pi^Y_{Y\cap
Z}$ and $\pi^Z_{Y\cap Z}$ in the following lemma.

\bn
\begin{lms}\label{aaa} (a) There
are isomorphisms of $\Ch(X)$-algebras
    \[B\oplus_{\pi,\gamma\pi} E \cong B\oplus_{\pi\eta,\pi\psi}E \cong
     \eta(B)\oplus_{Y\cap Z} \psi (E),\]
    where the second isomorphism is given by the map
    $\chi:B\oplus_{\pi\eta,\pi\psi}E \to A$  induced by the pair
    $(\eta,\psi)$.  Its components $\chi_x$ can be identified with $\psi_x$
    for $x\in Z$ and
    with $\eta_x$ for $x\in X\setminus Z$.

(b) Condition \eqref{inclusion} is equivalent to $\psi(E)\subset
    \pi_Z\big(A\oplus_Y \eta(B)\big)$.

(c) If $\fset$ is a finite subset of $A$ such that $\pi_Y(\fset)\subset_\ep
\eta(B)$ and $\pi_Z(\fset)\subset_\ep \psi(E)$, then $\fset\subset_\ep
\eta(B)\oplus_{Y\cap Z} \psi(E)=\chi\big(B\oplus_{\pi\eta,\pi\psi} E \big)
    $.
\end{lms}
\en

\bn
\begin{nots}
\label{d-bimultiplicative}
Let $\varphi:A \to B$ be a completely positive contractive (c.p.c.)  map between $C^{*}$-algebras; let $\Fh$ be a subset of $A$ and $\delta>0$. We say $\varphi$ is $\delta$-multiplicative on $\Fh$, or $(\Fh,\delta)$-multiplicative,  if $\|\varphi(ab)-\varphi(a)\varphi(b)\|<\delta$ for $a,b \in \Fh$. We say $\varphi$ is $\delta$-bimultiplicative on $\Fh$, if $\|\varphi(abc)-\varphi(a)\varphi(b)\varphi(c)\|<\delta$ for $a,b,c \in \Fh$.
\end{nots}
\en

\bn
\label{prop;approx-D-au}
\begin{props}
Let $\Dh$ be a separable unital strongly self-absorbing $C^{*}$-algebra. For any finite set $\mathcal{F}\subset \Dh$ and $\varepsilon>0$
there are a finite set $\mathcal{G}\subset \Dh$ and $\delta>0$
with the following property. For any unital $C^{*}$-algebra $A\cong A\otimes \Dh$ and any two
$(\mathcal{G},\delta)$-multiplicative u.c.p. maps $\varphi,\psi:\Dh \to A$
there is a unitary $u\in U(A)$ such that $\|\varphi(d)-u\psi(d)u^*\|<\varepsilon$ for all $d \in \mathcal{F}$. If $\Dh$ is additionally assumed to be $K_{1}$-injective, we may choose  $u\in U_{0}(A)$.
\end{props}
\en

\begin{nproof} Seeking a contradiction let us assume that for some given $\fset$ and $\ep$
there are no $\mathcal{G}$ and $\delta$ satisfying the conclusion of the Proposition. Let $(\mathcal{G}_n)$ be a sequence of increasing finite
subsets of $\Dh$ whose union is dense in $\Dh$ and let $\delta_n=1/n$ . Then there exist a sequence of $\Dh$-stable
unital $C^{*}$-algebras $(A_n)$ and two sequences  consisting of $(\mathcal{G}_n,\delta_n)$-multiplicative
u.c.p. maps $(\varphi_n)$
and $(\psi_n)$ with $\varphi_n,\psi_n:\Dh\to A_n$ such that for any $n$ and any $u_n\in U(A_n)$, $\|\varphi_n(d)-u_n\psi_n(d)u_n^*\|\geq\varepsilon$ for some $d \in \mathcal{F}_n$.

Set $B_n=\prod_{k\geq n} A_k$
and let $\nu_n:B_n\to B_{n+1}$ be the natural projection. Let us define
$\Phi_n,\Psi_n:\Dh\to B_n$ by
\[
\Phi_n(d)=(\varphi_n(d),\varphi_{n+1}(d),\dots)
\]
and
\[
\Psi_n(d)=(\psi_n(d),\psi_{n+1}(d),\dots),
\]
for all $d$ in $\Dh$. One
verifies immediately that if we set $\Phi=\varinjlim\,\Phi_n$, $\Psi=\varinjlim\,\Psi_n$ and $B=
\varinjlim\,(B_n,\nu_n)$, then $\Phi$ and $\Psi$ are unital $*$-homomorphisms from $\Dh$ to $B$. Since $\Dh$ is nuclear,
$B\otimes \Dh=
\varinjlim\,(B_n\otimes \Dh,\nu_n\otimes \id_{\Dh})$.
By applying \cite[Cor.~1.12]{TomsWinter:ssa} to $\Phi\otimes 1_{\Dh},\Psi\otimes 1_{\Dh}:\Dh \to B\otimes \Dh$, we find a unitary $V \in B \otimes \Dh$ such that  $\|\Phi(d)\otimes 1_{\Dh}-V \,(\Psi(d)\otimes 1_{\Dh})\,V^*\|<\varepsilon/4$ for all $d\in \mathcal{F}$.  One then  finds a sequence of unitaries
$V_n \in B_n\otimes \Dh$ such that $\|\Phi_n(d)\otimes 1_{\Dh}-V_n \,(\Psi_n(d)\otimes 1_{\Dh})\,V_n^*\|<\varepsilon/2$ for all $d\in \mathcal{F}$ and sufficiently large $n$. Projecting to  $A_n\otimes \Dh$,  we find that if $v_n\in U(A_n\otimes D)$
denotes the component of $V_n$ in $A_n\otimes \Dh$, then
\begin{equation}\label{1111}
\|\varphi_n(d)\otimes 1_{\Dh}-v_n \,(\psi_n(d)\otimes 1_{\Dh})\,v_n^*\|<\varepsilon/2
\end{equation}
 for all $d\in \mathcal{F}$.
By \cite[Prop.~1.9]{TomsWinter:ssa} there
is a sequence of unital $*$-homomorphisms $\theta_n: \Dh\ot \Dh\to \Dh$ such that
$\lim_{n\to \infty}\|\theta_n(d\ot 1_{\Dh})-d\|=0$ for all $d\in \Dh$.
Then $\gamma_n:A\ot \Dh \ot \Dh\to A\ot \Dh$, $\gamma_n=\id_A\ot \theta_n$
is a sequence of unital $*$-homomorphisms such that
$\lim_{n\to \infty}\|\gamma(a\ot 1_{\Dh})-a\|=0$ for all $a\in A\ot \Dh$.
Therefore, since $A_n$ is $\Dh$-stable, there
is a unital $*$-homomorphism $\gamma:A_n\otimes \Dh \to A_n$ such that $\|\gamma(a\ot
1_{\Dh})-a\|<\varepsilon/4$ for all $a$ of the form $\varphi_n(d)$ and $\psi_n(d)$ with $d\in \mathcal{F}$.
From ~\eqref{1111} we see that
             $\|\gamma(\varphi_n(d)\otimes 1_{\Dh})
             -\gamma(v_n) \,\gamma(\psi_n(d)\otimes
1_{\Dh})\,\gamma(v_n)^*\|<\varepsilon/2$   and hence        $\|\varphi_n(d)-\gamma(v_n)\psi_n(d)\gamma(v_n)^*\|<\varepsilon$
for all $d\in \mathcal{F}$. This contradicts our initial assumption.

In the $K_{1}$-injective case, to reach a contradiction we assume that the $u_{n}$ above are in $U_{0}(A_{n})$. By Theorem \ref{D-asu}, we may then assume that $V \in U_{0}(B \otimes \Dh)$; it is straightforward to show that the $V_{n}$ may also be chosen in $U_{0}(B_{n} \otimes \Dh)$. But then, $v_{n}$ and $\gamma(v_{n})$ are connected to the respective identities as well.
\end{nproof}

\bn
\begin{lms}\label{Lemma:Choi-Effros}
   Let $\Dh$ be a separable unital nuclear $C^{*}$-algebra.
For any finite set $\mathcal{G}_0\subset \Dh$ and $\delta_0>0$
there are a finite set $\mathcal{G}\subset \Dh$ and $\delta>0$
with the following property.
For any pair of unital $C^{*}$-algebras $1_B\in A\subset B$,
and any $(\mathcal{G},\delta)$-multiplicative u.c.p. map $\varphi:\Dh \to B$ satisfying
$\varphi(\mathcal{G})\subset_{\delta} A$, there is
a $(\mathcal{G}_0,\delta_0)$-multiplicative u.c.p. map  $\mu:\Dh \to A$ such that
$\|\varphi(d)-\mu(d)\|<\delta_0$ for any $d \in \mathcal{G}_0$.
\end{lms}

\begin{nproof} This is a known consequence of the Choi-Effros lifting theorem. One  proves this by contradiction along the lines of the
proof of the previous proposition.
\end{nproof}
\en

\bn
\begin{props}\label{prop:wsp-for-ssa}
Let $\Dh$ be a separable unital strongly self-absorbing $C^{*}$-algebra.
For any finite set $\mathcal{F}\subset \Dh$ and $\varepsilon>0$
there are a finite  set $\mathcal{G}\subset \Dh$ and $\delta>0$
with the following property.
For any pair of unital $C^{*}$-algebras $1_B\in A\subset B$ with $A\cong A\otimes \Dh$,
and any $(\mathcal{G},\delta)$-multiplicative u.c.p. map $\varphi:\Dh \to B$ satisfying
$\varphi(\mathcal{G})\subset_{\delta} A$, there is
a unital $*$-homomorphism $\psi:\Dh \to A$ such that $\|\varphi(d)-\psi(d)\|<\varepsilon$ for any $d \in \mathcal{F}$.
\end{props}

\begin{nproof} Let $\mathcal{G}_0$ and $\delta_0$ be given by Proposition~\ref{prop;approx-D-au} applied to $\Dh$, $\mathcal{F}$
and $\varepsilon/2$.
We may assume that $\mathcal{F}\subset \gset_0$ and $\delta_0<\ep/2$.  Let $\gset$ and $\delta$ be given by applying Lemma~\ref{Lemma:Choi-Effros} for $\Dh$, $\mathcal{G}_0$ and $\delta_0$. Suppose now that $\varphi$ is as in the statement and let $\mu:\Dh \to A$ be the perturbation of $\varphi$ yielded by Lemma~\ref{Lemma:Choi-Effros}.
Thus $\mu$ is $(\mathcal{G}_0,\delta_0)$-multiplicative and
$\|\varphi(d)-\mu(d)\|<\delta_0<\ep/2$ for all $d\in \fset\subset \gset_0$.
Since $A\cong A\otimes \Dh$, there is a  unital $*$-homomorphism $\eta:\Dh\to A$. Moreover,
by Proposition~\ref{prop;approx-D-au}, there is a unitary $u\in A$
such that $\|\mu(d)-u \,\eta(d)\,u^*\|<\varepsilon/2$ for all $d\in
\mathcal{F}$.
Let us set $\psi=u \,\eta\,u^*:\Dh \to A.$ Then
\[\|\varphi(d)-\psi(d)\|\leq \|\varphi(d)-\mu(d)\|+\|\mu(d)-u \,\eta(d)\,u^*\|<\varepsilon/2+\varepsilon/2=\varepsilon\]
for all $d\in \mathcal{F}$.
\end{nproof}
\en

 Proposition~\ref{prop:wsp-for-ssa} shows that separable unital strongly self-absorbing  $C^{*}$-algebras satisfy a relative
 perturbation property similar to weak semiprojectivity.
This is one of the key tools used in the proof of Theorem~\ref{thm:ssa-fields-are-trivial}.

\bn
\begin{lms}\label{lemma:local-D-sections}
Let $\Dh$ be a separable unital strongly self-absorbing  $C^{*}$-algebra. Let $X$ be a compact metrizable space and let $A$ be a unital $\Dh$-stable
   $\Ch(X)$-algebra with all fibres isomorphic to $\Dh$. Let $\mathcal{F}\subset A$ be
a finite set and let $\varepsilon>0$. For any $x\in X$
 there exist a closed neighborhood $U$ of $x$  and a unital $*$-homomorphism $\varphi: \Dh
\to A(U)$  such that $\pi_U(\mathcal{F})\subset_\varepsilon \varphi
(\Dh)$.
\end{lms}

\begin{nproof}
Since $A$ is unital and $\Dh$-stable there is a unital
$*$-homomorphism $\psi:\Dh \to A$.
By \cite[Cor.~1.12]{TomsWinter:ssa}, for each $x\in X$, the map $\psi_x:\Dh \to A(x)$ is approximately unitarily equivalent to some $*$-isomorphism
$\Dh\to A(x)$. Therefore there is a unitary $v\in A(x)$ such that $\pi_x(\mathcal{F})\subset_\varepsilon v \,\psi_x(\Dh)\,v^*$. Using the semiprojectivity of $\Ch(\mathbb{T})$ we lift  $v$ to a  unitary in $u\in A(V)$ for some closed neighborhood $V$ of $x$ and so
$\psi_x(\mathcal{F})\subset_\varepsilon \pi_x(u) \,\pi_x\psi(\Dh)\,\pi_x(u)^*$. Using also the upper continuity
of the norm in $\Ch(X)$-algebras,
 we conclude that there is closed neighborhood $U\subset V$ of $x$  such that $\pi_U(\mathcal{F})\subset_{\varepsilon} \pi_U(u)\,\pi_U\psi(\Dh)\, \pi_U(u)^*$.
\end{nproof}
\en

The following lemma is useful for  constructing  fibered morphisms,
see ~\ref{fibred-morphisms}.

\bn
\begin{lms}\label{lemma:multi-perturbation}
Let $\Dh$ be a separable unital strongly self-absorbing $K_1$-injective $C^{*}$-algebra.
 Let $(D_j)_{j\in J}$ be a finite family of $C^{*}$-algebras isomorphic to $\Dh$.
  Let $\varepsilon>0$ and for each $j\in J$ let
$\mathcal{H}_j\subset D_j$ be a finite set.
 Let $\mathcal{G}_j\subset D_j$ and
$\delta_j>0$ be given by Proposition~\ref{prop:wsp-for-ssa} applied to
$D_j$, $\mathcal{H}_j$ and $\varepsilon/2$.
    Let $X$ be
a compact metrizable space, let $(Z_j)_{j\in J}$ be  disjoint nonempty closed
subsets of $X$ and let $Y$ be a closed nonempty subset of $X$ such that
$X=Y\cup \big(\cup_j Z_j\big)$. Let $A=A(X)$ be a unital separable $\Dh$-stable $\Ch(X)$-algebra and let $B(Y)$ be a be a unital separable $\Dh$-stable $\Ch(Y)$-algebra.
Let $\mathcal{F}$ be a finite subset of $A$. Let $\eta:B(Y)\rightarrow A(Y)$ be a $\Ch(Y)$-linear unital
 $*$-monomorphism.  Suppose that $\varphi_j:D_j \to A(Z_j)$, $j\in J$,
 are unital $*$-homomorphisms  which satisfy the following conditions:
\begin{itemize}
        \item[(i)] $\pi_{Z_j}(\mathcal{F})\subset_{\varepsilon/2} \varphi_j(\mathcal{H}_j)$, for all $j \in
        J$,
        \item[(ii)] $\pi_Y(\mathcal{F})\subset_{\varepsilon} \eta(B)$,
        \item[(iii)] $\pi^{Z_j}_{Y\cap Z_j}\varphi_j(\mathcal{G}_j)
        \subset_{\delta_j} \pi^Y_{Y\cap Z_j}\eta(B)$, for all $j \in J$.
    \end{itemize}
Then, there are $\Ch(Z_j)$-linear unital $*$-monomorphisms $\psi_j:\Ch(Z_j)\ot D_j\to A(Z_j)$,
satisfying \begin{equation}\label{perturb-on-H}
    \|\varphi_j(c)-\psi_j(c)\|<\varepsilon/2, \,\,\text{for all}\,\,
    c\in \mathcal{H}_j,\,\,\text{and}\,\, j\in J,
\end{equation}
and such that if we set $E=\bigoplus_j \Ch(Z_j)\otimes D_j$, $Z=\cup_j Z_j$, and
$\psi:E \to A(Z)=\bigoplus_j A(Z_j)$, $\psi=\oplus_j\psi_j$, then $\pi^Z_{Y\cap
Z} (\psi(E))\subseteq \pi^Y_{Y\cap Z} (\eta(B))$, $\pi_Z(\mathcal{F})\subset_\varepsilon
\psi(E)$ and hence
\[\mathcal{F} \subset_{\varepsilon}\eta(B)\oplus_{Y\cap Z} \psi(E)=\chi( B
\oplus_{\pi\eta,\pi\psi} E),
 \]
where $\chi$ is the  isomorphism   induced by the pair $(\eta,\psi)$.
Moreover $B\oplus_{\pi\eta,\pi\psi} E$ is $\Dh$-stable.
\end{lms}

\begin{nproof} Let $\mathcal{F}=\{a_1,\dots,a_r\}\subset A$ be as in the statement. By
 (i), for each $j\in J$ we find
 $\{c_1^{(j)},\dots,c_r^{(j)}\}\subseteq \mathcal{H}_j$ such that
$\|\varphi_j(c_i^{(j)})-\pi_{Z_j}(a_i)\|<\varepsilon/2$ for all $i$.  The
$\Ch(X)$-algebra $A\oplus_Y\eta (B)\subset A$
is an extension of separable $\Dh$-stable $C^{*}$-algebras by \eqref{use}, and hence it
 is $\Dh$-stable by \cite[Thm.~4.3]{TomsWinter:ssa}.
 In particular  $\pi_{Z_j}(A\oplus_Y\eta (B))$ is $\Dh$-stable for each $j\in J$.

By (iii) and
Lemmas~\ref{a(x)delta}, ~\ref{sumover-Y}
 we obtain
\begin{equation*}\label{101iii} \varphi_j(\mathcal{G}_j)\subset_{\delta_j}
\pi_{Z_j}(A\oplus_Y\eta (B)).
\end{equation*}
 Applying Proposition~\ref{prop:wsp-for-ssa} we perturb
$\varphi_j$ to a  $*$-homomorphism $\psi_j:D_j \to \pi_{Z_j}(A\oplus_Y\eta (B))$
satisfying ~\eqref{perturb-on-H}, and hence such that $\|
\varphi_j(c_i^{(j)})-\psi_j(c_i^{(j)})\|<\varepsilon/2,$  for all $i,j$. Therefore
\begin{equation*}\label{5}
\| \psi_j(c_i^{(j)})-\pi_{Z_j}(a_i)\|\leq \|
\psi_j(c_i^{(j)})-\varphi_j(c_i^{(j)})\|+\|\varphi_j(c_i^{(j)})-\pi_{Z_j}(a_i)\|<\ep.
\end{equation*}  This
shows that $\pi_{Z_j}(\fset)\subset_{\ep}\psi_j(D_j)$.
 We extend $\psi_j$ to a
$\Ch(Z_j)$-linear  $*$-monomorphism $\psi_j:\Ch(Z_j)\otimes D_j \to
 \pi_{Z_j}(A\oplus_Y\eta (B))$ and then we define $E$, $\psi$ and $Z$
  as in the statement. In this way
we obtain that $\psi:E \to (A\oplus_Y\eta(B))(Z)\subset A(Z)$ satisfies
\begin{equation}\label{abc}
   \pi_Z(\mathcal{F})\subset_{\varepsilon} \psi(E).
\end{equation}
The property $\psi(E)\subset (A\oplus_Y\eta(B))(Z)$ is equivalent to
$\pi^Z_{Y\cap Z} (\psi(E))\subset \pi^Y_{Y\cap Z} (\eta(B))$ by
Lemma~\ref{aaa}(b). Finally, from (ii), \eqref{abc} and Lemma~\ref{aaa}\,(c) we
get $\mathcal{F} \subset_{\varepsilon} \eta(B)\oplus_{Y\cap Z}
  \psi(E).$
  One verifies the $\Dh$-stability of $B\oplus_{\pi\eta,\pi\psi} E$
by repeating the argument that showed the $\Dh$-stability of $A\oplus_Y\eta (B)$.
\end{nproof}
\en

\bn\label{fibred-morphisms}
Let $\Dh$ be a separable unital strongly self-absorbing $K_1$-injective $C^{*}$-algebra. Let $A$ be a
    $\Ch(X)$-algebra,
  let $\fset\subset A$ be a finite set and let $\ep>0$.
    An $(\fset,\ep,\Dh)$-\emph{approximation} of $A$ is a finite collection
\begin{equation}\label{FE-approx}
\alpha=\{\fset,\ep,\{U_i,\varphi_i:D_i\to
A(U_i),\mathcal{H}_i,\mathcal{G}_i,\delta_i\}_{i\in I}\},
\end{equation}
 with the following properties: $(U_i)_{i\in I}$ is a finite
family of closed subsets of $X$  whose interiors  cover  $X$;
 $(D_i)_{i\in I}$ are $C^{*}$-algebras isomorphic to  $\Dh$;
 for each $i\in I$,
$\varphi_i:D_i\to A(U_i)$ is a unital $*$-homomorphism;
  $\mathcal{H}_i\subset D_i$ is a finite
  set such that $\pi_{U_i}(\fset)\subset_{\ep/2} \varphi_i(\mathcal{H}_i)$; the finite set
$\mathcal{G}_i\subset D_i$ and $\delta_i>0$ are given by
Proposition~\ref{prop:wsp-for-ssa} applied to  $D_i$ for the input data $\mathcal{H}_i$ and $\ep/2$.

 It is
useful to consider the following operation of restriction. Suppose that $Y$ is
a closed subspace of $X$ and let $(V_j)_{j\in J}$ be a finite family of closed
subsets of $Y$ which refines the family $(Y\cap U_i)_{i\in I}$ and such that
the interiors of the $V_j\,'s$ form a cover of $Y$. Let $\iota:J\to I$ be a map
such that $V_j\subseteq Y\cap U_{\iota(j)}$. Define
\[\iota^*(\alpha)=\{\pi_Y(\fset),\ep,\{V_j,\pi_{V_j}\varphi_{\iota(j)}:D_{\iota(j)}\to
A(V_j),\mathcal{H}_{\iota(j)},\mathcal{G}_{\iota(j)},\delta_{\iota(j)}\}_{j\in
J}\}.\] It is obvious that  $\iota^*(\alpha)$ is a
$(\pi_Y(\fset),\ep,\Dh)$-approximation of $A(Y)$. The operation
$\alpha\mapsto \iota^*(\alpha)$  is useful even  in the case $Y=X$. Indeed, by
applying this procedure we can refine the cover of $X$ that appears in a given
$(\fset,\ep,\Dh)$-approximation of $A$.

    An
$(\fset,\ep,\Dh)$-approximation $\alpha$ 
  is subordinate to an
$(\fset',\ep',\Dh)$-approximation,
$\alpha'=\{\fset',\ep',\{U_{i'},\varphi_{i'}:D_{i'}\to
A(U_{i'}),\mathcal{H}_{i'},\mathcal{G}_{i'},\delta_{i'}\}_{{i'}\in I'}\},$
written $\alpha\prec\alpha'$, if
\begin{itemize}
        \item[(i)] $\fset\subseteq \fset'$,
        \item[(ii)] $\varphi_i(\mathcal{G}_i)\subseteq \pi_{U_i}(\fset')$
        for all $i\in I$, and
        \item[(iii)] $\ep'< \min\big(\{\ep\}\cup\{\delta_i,\, i\in I\}\big)$.
    \end{itemize}
Let us note that, with notation as above,  we have $\iota^*(\alpha)\prec
\iota^*(\alpha')$ whenever $\alpha\prec\alpha'$.
\en

\bn \label{fibred-monomorphism}
Let us recall some terminology from \cite{Dad:bundles-fdspaces}.
 A  $\Ch(Z)$-algebra $E$ is called
  $\Dh$-\emph{elementary} if $E\cong \Ch(Z)\ot \Dh$.
  Let $A$ be a unital $\Ch(X)$-algebra.
 A unital  $n$-\emph{fibered $\Dh$-monomorphism
$(\psi_0,\dots,\psi_n)$ into $A$}
      consists of $(n+1)$ unital $*$-monomorphisms of $\Ch(X)$-algebras
$\psi_i:E_i \to A(Y_i)$, where $Y_0,\dots,Y_n$ is a closed cover of $X$,
    each $E_i$ is a
   $\Dh$-elementary  $\Ch(Y_i)$-algebra and
    \begin{equation*}\label{containement-fibred}
\pi^{Y_i}_{Y_i\cap
    Y_j}\psi_i(E_i)\subseteq \pi^{Y_j}_{Y_i\cap
    Y_j}\psi_j(E_j), \quad \text{ for all}\,\, 0\leq i\leq j \leq n.
\end{equation*}
    Given an $n$-fibered morphism into $A$ we have an associated \emph{continuous} $\Ch(X)$-algebra
    defined  as the  fibered product (pullback)  of the $*$-monomorphisms $\psi_i$:
    \begin{equation*}
    A(\psi_0,\dots,\psi_n)=\{(d_0,\dots d_n)\,:\,\, d_i\in E_i,\,
    \pi^{Y_i}_{Y_i\cap Y_j}\psi_i(d_i)=\pi^{Y_j}_{Y_i\cap Y_j}\psi_j(d_j)
    \,\,\text{for all}\, i,j\}
\end{equation*}
and an induced $\Ch(X)$-linear monomorphism
$$\eta=\eta_{(\psi_0,\dots,\psi_n)}:A(\psi_0,\dots,\psi_n)\to A\subset
\bigoplus_{i=0}^n A(Y_i), $$
\[\eta(d_0,\dots d_n)=\big(\psi_0(d_0),\dots,\psi_n(d_n)\big).\]

 Let us set
$X_k=Y_k\cup\cdots\cup Y_n$. Then $(\psi_k,\dots\psi_n)$ is an $(n-k)$-fibered
$\Dh$-monomorphism into $A(X_k)$. Let $\eta_k :
A(X_k)(\psi_k,\dots\psi_n)\to A(X_k)$ be the induced map and set
$B_k=A(X_k)(\psi_k,\dots\psi_n)$. Let us note that $B_0=A(\psi_0,\dots,\psi_n)$
and that there are natural $\Ch(X_{k-1})$-isomorphisms
\begin{equation}\label{eq-recurrence-B(k)}
B_{k-1}\cong B_k\oplus_{\pi\eta_k,\pi\psi_{k-1}} E_{k-1}\cong
B_k\oplus_{\pi,\gamma_k\pi} E_{k-1},
\end{equation}
where $\pi$ stands for $\pi_{X_k\cap Y_{k-1}}$ and $\gamma_k: E_{k-1}(X_k\cap
Y_{k-1}) \to B_k(X_k\cap Y_{k-1})$ is defined by
$(\gamma_k)_x=(\eta_k)_x^{-1}(\psi_{k-1})_x$, for all $x\in X_k\cap Y_{k-1}$.

We say that a unital $n$-fibered $\Dh$-monomorphism
$(\psi_0,\dots,\psi_n)$ into $A$ is \emph{locally extendable} if
 there is another unital $n$-fibered
$\Dh$-monomorphism $(\psi'_0,...,\psi'_n)$ into $A$ such that
$\psi'_k:E_k^\prime\to A(Y_k')$,  $Y_k'$ is a closed neighborhood of $Y_k$,
$E_k'(Y_k)=E_k$ and
$\pi_{Y_k}\psi'_k=\psi_k$, $k=0,...,n.$

\en
The following proposition is a crucial technical result in this section.

\bn
\begin{props}\label{basic-approx}
Let $\Dh$ be a separable unital strongly self-absorbing $K_1$-injective $C^{*}$-algebra.
  Let $X$ be a finite dimensional compact metrizable space
    and let $A$ be a unital separable continuous
    $\Ch(X)$-algebra the fibres of which are isomorphic to $\Dh$.
 For any finite set $\fset\subset A$
 and any $\ep>0$ there exist $n\leq\mathrm{dim}(X)$ and an $n$-fibered unital $\Dh$-monomorphism
$(\psi_0,\dots,\psi_n)$ into $A$   which induces a unital $*$-monomorphism
$\eta:A(\psi_0,...,\psi_n)\to A$ such that  $\fset\subset_\ep \eta
(A(\psi_0,...,\psi_n))$.
\end{props}

\begin{nproof} The $C^{*}$-algebra $A$ is $\Dh$-stable by \cite{HirshbergRordamWinter:absorb-ssa}.
By Lemma~\ref{lemma:local-D-sections}, Proposition~\ref{prop:wsp-for-ssa} and the compactness of $X$, for any finite set $\fset\subset A$
 and any $\ep>0$
  there is an $(\fset,\ep,\Dh)$-approximation of $A$. Moreover, for
  any finite set $\fset\subset A$,
  any $\ep>0$ and any $n$, there is a sequence
  $\{\alpha_k:\, 0\leq k\leq n\}$ of $(\fset_k,\ep_k,\Dh)$-approximations
  of $A$ such that $(\fset_0,\ep_0)=(\fset,\ep)$ and $\alpha_k$ is
   subordinated to $\alpha_{k+1}$:
   \[\alpha_0\prec\alpha_1\prec\cdots\prec\alpha_n.\]
Indeed, assume that $\alpha_k$ was constructed. Let us choose a finite set
$\fset_{k+1}$ which contains $\fset_k$ and liftings to $A$ of all the elements
in $\bigcup_{i_k\in I_k}\varphi_{i_k}(\gset_{i_k}).$ This choice takes care of
the above conditions (i) and (ii). Next we choose $\ep_{k+1}$ sufficiently
small such that (iii) is satisfied. Let $\alpha_{k+1}$ be an
$(\fset_{k+1},\ep_{k+1},\Dh)$-approximation  of $A$ given by
Lemma~\ref{lemma:local-D-sections} and Proposition~\ref{prop:wsp-for-ssa}. Then obviously $\alpha_k\prec\alpha_{k+1}$. Fix
a tower of approximations of $A$ as above where $n=\mathrm{dim}(X)$.

We may assume that the set $X$ is infinite.
By \cite[Lemma~3.2]{BK:bundles}, for every open cover $\mathcal{V}$ of $X$
there is a finite open cover $\mathcal{U}$ which refines $\mathcal{V}$ and such
that the set $\mathcal{U}$ can be partitioned into $n+1$ nonempty subsets
consisting of elements with pairwise disjoint closures.
 Since we can refine simultaneously
the covers that appear in a finite family $\{\alpha_k:\, 0\leq k\leq n\}$ of
approximations while preserving subordination,
 we may arrange not only that all $\alpha_k$ share the same cover $(U_i)_{\in I},$
but moreover, that the cover $(U_i)_{i\in I}$ can be partitioned into $n+1$
subsets $\U_0,\dots,\U_n$ consisting of mutually disjoint elements. For
definiteness, let us write $\U_k=\{U_{i_k}:i_k\in I_k\}$. Now for each $k$ we
consider the closed subset of $X$
\[Y_k=\bigcup_{i_k\in I_k}\,U_{i_k},\]
the map $\iota_k:I_k\to I$ and the
$(\pi_{Y_{k}}(\fset_k),\ep_k,\Dh)$-approximation of $A(Y_k)$ induced
by $\alpha_k$,  which is of the form
\[ \iota_k^*(\alpha_k)=\{\pi_{Y_k}(\fset_k),\ep,\{U_{i_k},\varphi_{i_k}:D_{i_k}\to
A(U_{i_k}),\mathcal{H}_{i_k},\mathcal{G}_{i_k},\delta_{i_k}\}_{{i_k}\in
I_k}\},\] where each $U_{i_k}$ is nonempty.
 We have
\begin{equation}\label{x0}
   \pi_{U_{i_k}}(\fset_k)\subset_{\ep_k/2}\varphi_{i_k}(\mathcal{H}_{i_k}),
\end{equation}
 by construction. Since $\alpha_k\prec\alpha_{k+1}$ we also have
\begin{equation}\label{x1}
    \fset_k\subseteq \fset_{k+1},
\end{equation}
\begin{equation}\label{x2}
    \varphi_{i_k}(\gset_{i_k})\subseteq\pi_{U_{i_k}}(\fset_{k+1}),\,\,\text{for
    all}\,\, i_k\in I_k,
\end{equation}
\begin{equation}\label{x3}
   \ep_{k+1}< \min\big(\{\ep_k\}\cup\{\delta_{i_k},\, {i_k}\in I_k\}\big).
\end{equation}
Set $X_k=Y_k\cup \dots \cup Y_n$ and $E_k=\oplus_{i_k} \,\Ch(U_{i_k})\ot D_{i_k}$
for $0\leq k \leq n$.
 We shall construct  a sequence  of unital $\Ch(Y_k)$-linear $*$-monomorphisms, $\psi_k:E_k \to
A(Y_k)$, $k=n,...,0$,  such that $(\psi_k,\dots,\psi_n)$ is an $(n-k)$-fibered
monomorphism into $A(X_k)$. Each map
\[\psi_k=\oplus_{i_k}\psi_{i_k}:E_k \to A(Y_k)=\oplus_{i_k}\,A(U_{i_k})\]
will have components
 $\psi_{i_k}:\Ch(U_{i_k})\ot D_{i_k}\to
A(U_{i_k})$  whose restrictions to $D_{i_k}$ will be perturbations of
$\varphi_{i_k}:D_{i_k}\to A(U_{i_k})$, ${i_k}\in I_k$.
 We shall construct the maps $\psi_k$  by induction on decreasing $k$
  such that if
 $ B_{k}=A(X_k)(\psi_k,\dots,\psi_n)$ and  $\eta_{k}:B_{k}\to A(X_k)$
 is the map induced by the
 $(n-k)$-fibered monomorphism $(\psi_k,\dots,\psi_n)$, then
\begin{equation}\label{6k}
    \pi_{X_{k+1}\cap U_{i_k}}\big(\psi_{i_k}(D_{i_k})\big)
    \subset \pi_{X_{k+1}\cap U_{i_k}}
    \big(\eta_{k+1}(B_{k+1})\big), \,\forall \,{i_k} \in I_k,
\end{equation} and
\begin{equation}\label{7k}
    \pi_{X_k}(\fset_k)\subset_{\ep_k}\eta_{k}(B_{k}).
\end{equation}
Note that ~\eqref{6k} is equivalent to
\begin{equation}\label{6kk}
    \pi_{X_{k+1}\cap Y_k}\big(\psi_k(E_k)\big)\subset \pi_{X_{k+1}\cap Y_k}
    \big(\eta_{k+1}(B_{k+1})\big).
\end{equation}
Let us note that $B_k$ is $\Dh$-stable by \cite {HirshbergRordamWinter:absorb-ssa}, since each of its fibres are $\Dh$-stable and $X$ is finite dimensional.

 For the first step of
induction, $k=n$, we choose  $\psi_n=\oplus_{i_n}\widetilde{\varphi}_{i_n}$
where $\widetilde{\varphi}_{i_n}:\Ch(U_{i_n})\ot D_{i_n}\to A(U_{i_n})$ are
$\Ch(U_{i_n})$-linear extensions of the original $\varphi_{i_n}$. Then we set $B_n= E_n$
and $\eta_n=\psi_n$.
 Assume  now that
$\psi_n,\dots,\psi_{k+1}$ were constructed. Next we  shall construct  $\psi_{k}$. Condition~\eqref{7k} formulated
for $k+1$ becomes
\begin{equation}\label{7kk}
    \pi_{X_{k+1}}(\fset_{k+1})\subset_{\ep_{k+1}}\eta_{k+1}(B_{k+1}).
\end{equation}
  Since $\ep_{k+1}<\delta_{i_k}$, by using  \eqref{x2} and
~\eqref{7kk}  we obtain
\begin{equation}\label{iii}
\pi_{X_{k+1}\cap
U_{i_k}}\big(\varphi_{i_k}(\gset_{i_k})\big)\subset_{\delta_{i_k}}
\pi_{X_{k+1}\cap U_{i_k}}
    \big(\eta_{k+1}(B_{k+1})\big), \,\,\text{for all}\, \,{i_k} \in I_k.
\end{equation}
  Since $\fset_k\subseteq \fset_{k+1}$
 and $\ep_{k+1}<{\ep_k}$, condition~\eqref{7kk} gives
\begin{equation}\label{ii}
    \pi_{X_{k+1}}(\fset_{k})\subset_{\ep_k}\eta_{k+1}(B_{k+1}).
\end{equation}
 Conditions \eqref{x0},
\eqref{iii}  and  \eqref{ii} enable us to apply
Lemma~\ref{lemma:multi-perturbation} and perturb $\varphi_{i_k}$ to  unital
$*$-homomorphisms $D_{i_k}\to A(U_{i_k})$ whose $\Ch(U_{i_k})$-linear extension $\psi_{i_k}:\Ch(U_{i_k})\ot D_{i_k}\to A(U_{i_k})$ satisfy
 ~\eqref{6k} and ~\eqref{7k}.
 We set
 $\psi_k=\oplus_{i_k}\,\psi_{i_k}$ and
  this completes the construction of $\Psi_k$ and hence of $(\psi_0,\dots,\psi_n)$.
  Condition ~\eqref{7k} for $k=0$ gives
$\fset\subset_\ep \eta_0(B_0)=\eta(A(\psi_0,\dots,\psi_n)).$ Thus
$(\psi_0,\dots,\psi_n)$ satisfies the conclusion of the theorem. Finally let us
note that it can happen that $X_k=X$ for some $k>0$. In this case
$\fset\subset_\ep A(\psi_k,...,\psi_n)$ and for this reason we write $n\leq
\mathrm{dim}(X)$ in the statement of the theorem.
\end{nproof}
\en

\bn
\begin{rems}\label{polys} We can arrange that the $n$-fibered morphism $(\psi_0,\dots,\psi_n)$ from the conclusion
of Proposition~\ref{basic-approx}
 is locally extendable. Fix a metric $d$ for the topology of $X$. Then we may arrange that
there is a closed cover $\{Y'_0,...,Y'_n\}$ of $X$
 and a number $\ell>0$ such that
 $\{x:d(x,Y'_i)\leq \ell\}\subset Y_i$ for $i=0,...,n$.
 Indeed, when we choose the finite closed cover $\mathcal{U}=(U_i)_{i\in I}$ of $X$
 in the proof of
 Proposition~\ref{basic-approx} which can be partitioned into $n+1$
subsets $\U_0,\dots,\U_n$ consisting of mutually disjoints elements,
  as given by \cite[Lemma~3.2]{BK:bundles},
 and which refines all the covers
  $\mathcal{U}(\alpha_0),...,\mathcal{U}(\alpha_n)$
  corresponding to $\alpha_0,...,\alpha_n$,
 we may assume that $\mathcal{U}$ also refines the covers given by the
 interiors of  the elements of
 $\mathcal{U}(\alpha_0),...,\mathcal{U}(\alpha_n)$.
 Since each $U_i$ is compact and $I$ is finite, there is $\ell>0$ such that
 if $V_i=\{x:d(x,U_i)\leq \ell\}$, then the cover $\mathcal{V}=(V_i)_{i\in I}$
 still refines all of $\mathcal{U}(\alpha_0),...,\mathcal{U}(\alpha_n)$
 and for each $k=0,...,n$, the elements
 of $\mathcal{V}_k=\{V_i:U_i\in \mathcal{U}_k\}$,
 are still mutually
 disjoint. We shall use the cover
 $\mathcal{V}$ rather than $\mathcal{U}$ in the proof of Proposition~\ref{basic-approx}
 and observe that $Y'_k\stackrel{def}{=}\bigcup_{i_k\in I_k} U_{i_k}\subset
 \bigcup_{i_k\in I_k} V_{i_k}=Y_k$ has the desired property.
 Finally let us note that if we define $\psi_i':E_i(Y_i')\to A(Y_i')$
 by $\psi_i'=\pi_{Y_i'}\psi_i$, then $(\psi_0',\dots,\psi_n')$
 is a locally extendable $n$-fibered unital $\Dh$-monomorphism into $A$ which satisfies
 the conclusion of Proposition~\ref{basic-approx} since $\pi_{Y'_i}(\fset)\subset_\ep \psi_i'(E_i(Y_i'))$
  for all $i=0,\dots,n$ and $X=\bigcup_{i=1}^{\,n}Y_i'$.
\end{rems}
\en

For a unital $C^{*}$-algebra $\Dh$ we let $\mathrm{End}(\Dh)$ denote the space of all
unital $*$-endomorphisms of $\Dh$ endowed with the point-norm topology. Throughout the remainder of
this section we assume that
 $X$  is a compact metrizable space and that $\Dh$ is a unital, separable,
 $K_{1}$-injective and
strongly self-absorbing $C^{*}$-algebra.

 \bn
\begin{lms}\label{lemma:pre-extend} Any continuous map $\alpha: X\times\{0\}\to \mathrm{End}(\Dh)$ extends to a
continuous map $\Phi:X\times [0,1]\to \mathrm{End}(\Dh)$ such that $\Phi_{(x,1)}=\id_{\Dh}$ and $\Phi_{(x,t)}\in
\Aut(\Dh)$ for all $x\in X$ and $t\in (0,1]$.
\end{lms}

\begin{nproof}
Let us identify $\alpha$ with a unital $*$-homomorphism $\alpha:\Dh\to \Ch(X)\ot \Dh$.
By Theorem~\ref{D-asu} there is a unitary-valued continuous map $(0,1]\to
U(\Ch(X)\otimes \Dh)$, $t\mapsto u_t$, with $u_1=1_{\Ch(X)\otimes \Dh}$ and such that
\begin{equation*}\label{d(Y)}
\lim_{t \to 0} \|u_t d u_t^*-\alpha(d)\|=0,\,\,\text{for all}\quad d \in
 D.
\end{equation*}
Therefore the equation \[\Phi_{(x,t)}(d)=
 \left\{
\begin{array}{ll}
        \alpha_x(d), & \hbox{if $t=0$,} \\
        u_t(x)\,d \,u_t(x)^*, &\hbox{if $t\in (0,1]$,} \\
\end{array}
\right.\] defines a continuous map $\Phi:X\times [0,1]\to \mathrm{End}(\Dh)$ which
extends $\alpha$, $\Phi_{(x,1)}=\id_{\Dh}$, and such that $\Phi(X\times (0,1])\subset \mathrm{Aut}(\Dh)$.
\end{nproof}
\en

\bn
\begin{props}\label{prop:extending}
   Let $Y$ be a closed subset of
    $X$ and let $V$ be a closed neighborhood of $Y$.
    Suppose that a   map $\gamma:Y \to \mathrm{End}(\Dh)$
   extends to a continuous map $\alpha:V \to \mathrm{End}(\Dh)$.
   Then there is a continuous extension $\eta:X\to\mathrm{End}(\Dh)$
   of $\gamma$
    such that
    $\eta(X\setminus Y)\subset \mathrm{Aut}(\Dh)$.
\end{props}

\begin{nproof}
First we prove the proposition in the case  when $V=X$.
 By Lemma~\ref{lemma:pre-extend},
there exists a continuous map $\Phi:X\times [0,1]\to \mathrm{End}(\Dh)$ which
extends $\alpha$, i.e. $\Phi_{(x,0)}=\alpha_x$ for $x\in X$, and such that
$\Phi(X\times (0,1])\subset \mathrm{Aut}(\Dh)$. Let $d$ be a metric for the
topology of $X$ such that $\mathrm{diam}(X)\leq 1$. The equation
$\eta_x=\Phi_{(x,d(x,Y))}$ defines a map on $X$ that satisfies the conclusion of
the proposition.

Suppose now that $V$ is a closed neighborhood of $Y$.
By Lemma~\ref{lemma:pre-extend}, there is a continuous map (homotopy) $\Phi:V\times [0,1]\to \End(\Dh)$
such that $\Phi_{(x,0)}=\id_\Dh$ and $\Phi_{(x,1)}=\alpha_x$
for all $x\in V$. Let $\widehat\Phi$ be the extension of $\Phi$
to $V\times [0,1]\cup X \times \{0\}$ obtained by setting $\widehat\Phi_{(x,0)}=\id_\Dh$ for $x\in X\setminus{V}$.
 Let us
define $\lambda:X \to [0,1]$ by
$\lambda(x)=d(x,X\setminus V)\big(d(x,X\setminus V)+d(x,Y)\big)^{-1}$ and $\widehat\alpha:X\to \End(\Dh)$ by
$\widehat\alpha_{x}=\widehat\Phi_{(x,\lambda(x))}$ and observe that
$\widehat\alpha$ extends $\gamma$ to $X$. The conclusion follows now from the first part of the proof.
\end{nproof}
\en

If $\varphi:\Dh\to \Ch(X)\otimes \Dh$ is a $*$-homomorphism, we denote by
$\widetilde{\varphi}:\Ch(X)\otimes \Dh \to \Ch(X)\otimes \Dh$ its $\Ch(X)$-linear
extension.

 \bn
\begin{lms}\label{local-triviality-CX-semiprojective}
 Let $Y,$ $Z$ be closed subsets of $X$ such
that $X=Y\cup Z$.
     Let $\gamma:\Dh \to \Ch(Y\cap Z)\otimes \Dh$ be a
   unital
    $*$-homomorphism.
Assume that there is a closed neighborhood $V$ of $Y\cap Z$ in $Y$ and a
   unital
    $*$-homomorphism $\alpha:\Dh \to \Ch(V)\otimes \Dh$
     such that
    $\alpha_x=\gamma_x$ for all $x\in Y\cap Z$.
    Then the pullback
     $\Ch(Y)\otimes \Dh\oplus_{\pi_{Y\cap Z},\widetilde{\gamma}\pi_{Y\cap Z}} \Ch(Z)\otimes \Dh$
     is isomorphic to
    $\Ch(X)\otimes \Dh$.
\end{lms}

\begin{nproof} By
Prop.~\ref{prop:extending} there is a unital $*$-homomorphism $\eta:\Dh\to
\Ch(Y)\otimes \Dh$ such that $\eta_x=\gamma_x$ for $x \in Y\cap Z$ and such that
$\eta_x \in \mathrm{Aut}(\Dh)$ for $x\in Y\setminus Z$. One checks immediately that
the pair $(\widetilde\eta,
 \mathrm{id}_{\Ch(Z)\otimes \Dh})$ defines a $\Ch(X)$-linear
isomorphism: \newline$\Ch(X)\otimes \Dh=\Ch(Y)\otimes \Dh\oplus_{\pi_{Y\cap Z},\pi_{Y\cap
Z}} \Ch(Z)\otimes \Dh\to \Ch(Y)\otimes \Dh\oplus_{\pi_{Y\cap
Z},\widetilde{\gamma}\pi_{Y\cap Z}} \Ch(Z)\otimes \Dh$.
Indeed if $P$ denotes the later pullback, there is commutative diagram
\[\xymatrix{
0 \ar[r] & \Ch_0(Y\setminus Z)\ot \Dh\ar[r]\ar[d]^{\widetilde\eta_{Y\setminus Z}}\ & \Ch(X)\ot \Dh\ar[r]
 \ar[d]^{(\widetilde\eta,
 \mathrm{id})}& \Ch(Z)\ot \Dh
\ar[r]\ar@{=}[d]&0\\
0 \ar[r] & \Ch_0(Y\setminus Z)\ot \Dh\ar[r] & P\ar[r] & \Ch(Z)\ot \Dh
\ar[r]&0
}\]
and hence $(\widetilde\eta,
 \mathrm{id}_{\Ch(Z)\otimes \Dh})$ must be bijective, since
 $\widetilde\eta_{Y\setminus Z}$ is so. \end{nproof}
\en

\bn
\begin{cors}\label{local-triviality}  Let $Y$, $Z$ and $Z'$ be
closed subsets of a compact metrizable space $X$ such that $Z'$ is a neighborhood
of $Z$ and $X=Y\cup Z$.
  Let $B$ be a  $\Ch(Y)$-algebra
    isomorphic to $\Ch(Y)\otimes \Dh$
    and let $E=\Ch(Z')\otimes \Dh.$
     Suppose that $\alpha: E(Y\cap Z')\to B(Y\cap Z')$  is
 a unital $*$-monomorphism of $\Ch(Y\cap Z')$-algebras. If
$\gamma=\alpha_{Y\cap Z}$,
   then $B(Y)\oplus_{\pi_{\SSS Y\cap Z},\gamma\,\pi_{\SSS Y\cap Z}}E(Z)$ is
isomorphic to $\Ch(X)\ot \Dh$.
\end{cors}

\begin{nproof} This follows from Lemma~\ref{local-triviality-CX-semiprojective} since $Y\cap Z'$ is a closed
neighborhood of $Y\cap Z$ in $Y$.
\end{nproof}\en

\bn
\begin{props}\label{prop:local-triviality}
Let $A$ be a unital $\Ch(X)$-algebra. If
$(\psi_0,\dots,\psi_n)$ is  a locally extendable, unital $n$-fibered
$\Dh$-monomorphism  into $A$, then $A(\psi_0,\dots,\psi_n)$ is isomorphic
to $\Ch(X)\otimes \Dh$.
\end{props}

\begin{nproof}  By assumption, there exists 
another unital $n$-fibered $\Dh$-monomorphism $(\psi'_0,...,\psi'_n)$ into $A$
such that $\psi'_k:\Ch(Y'_k)\ot \Dh\to A(Y_k')$, $Y_k'$ is a closed neighborhood of
$Y_k$, and $\pi_{Y_k}\psi'_k=\psi_k$, $k=0,...,n.$ Let $X_k$, $B_k$, $\eta_k$ and
$\gamma_k$ be as in ~\ref{fibred-monomorphism}. We need to show that
$B_0$ is  isomorphic to $\Ch(X)\ot \Dh$. To this purpose we prove by induction
on decreasing $k$ that the $\Ch(X_k)$-algebras $B_k$ are  trivial. Indeed
$B_n=\Ch(X_n)\ot \Dh$ and assuming that $B_k$ is  trivial, it follows by
Corollary~\ref{local-triviality} that $B_{k-1}$ is  trivial, since by
\eqref{eq-recurrence-B(k)}
\begin{equation*}\label{eq-recurrence-B(k)*}
B_{k-1}\cong B_k\oplus_{\pi\eta_k,\pi\psi_{k-1}} E_{k-1}\cong
B_k\oplus_{\pi,\gamma_k\pi} E_{k-1},\quad(\,\, \pi=\pi_{X_k\cap Y_{k-1}})
\end{equation*}
and $\gamma_k: E_{k-1}(X_k\cap Y_{k-1}) \to B_k(X_k\cap Y_{k-1})$,
$(\gamma_k)_x=(\eta_k)_x^{-1}(\psi_{k-1})_x$, extends to a $*$-monomorphism
$\alpha: E_{k-1}(X_k\cap Y'_{k-1}) \to B_k(X_k\cap Y'_{k-1})$,
$\alpha_x=(\eta_k)_x^{-1}(\psi'_{k-1})_x$. \end{nproof} \en

\bn
We are now ready to complete the first proof of our main result:

\begin{nnproof} (of Theorem~\ref{thm:ssa-fields-are-trivial}): Set $B=\Ch(X)\ot \Dh$. By Propositions ~\ref{prop:local-triviality},
 ~\ref{basic-approx} and Remark~\ref{polys} there is a sequence $(\theta_k)_{k=1}^\infty$
consisting of unital $*$-monomorphisms of $\Ch(X)$-algebras $\theta_k:B\to A$
such that the sequence $(\theta_k(B))_{k=1}^\infty$ exhausts $A$ in the sense
that for any finite subset $\gset$ of $A$ and any $\delta>0$, $\gset\subset_\delta \theta_k(B)$ for some $k\geq 1$.
Using Proposition ~\ref{prop:wsp-for-ssa},
 after passing to a
subsequence of $(\theta_k)_{k=1}^\infty$, we construct a sequence of finite
sets $\fset_k\subset B$ and a sequence of unital $\Ch(X)$-linear  $*$-monomorphisms $\mu_k:B\to
B$ such that
\begin{itemize}
        \item[(i)] $\|\theta_{k+1}\mu_k(a)-\theta_{k}(a)\|<2^{-k}$ for
        all $a \in \fset_k$ and all $k\geq 1$,
        \item[(ii)] $\mu_k(\fset_k)\subset \fset_{k+1}$ for all $k\geq 1$,
        \item[(iii)] $\bigcup_{j=k+1}^\infty\,
        \big(\mu_{j-1}\circ\cdots\circ\mu_{k}\big)^{-1}(\fset_j)$ is dense
        in $B$  and
        $\bigcup_{j=k}^\infty\,\theta_{j}(\fset_j)$ is
        dense in $A$ for all $k\geq 1$.
    \end{itemize}
    Arguing as in the proof of \cite[Prop.~2.3.2]{Ror:encyclopedia},
one verifies that $$\Delta_k(a)=\lim_{j \to
\infty}\theta_j\circ\big(\mu_{j-1}\circ\cdots\circ\mu_{k}\big) (a)$$ defines a
sequence of $*$-monomorphisms  $\Delta_k:B \to A$  such that
$\Delta_{k+1}\mu_k=\Delta_k$ and the induced map   $\Delta:\varinjlim_k
(B,\mu_k)\to A$ is an isomorphism of $\Ch(X)$-algebras. Let us show that
$\varinjlim_k (B,\mu_k)$ is isomorphic to $B$. To this purpose, in view of
Elliott's intertwining argument, it suffices to show that each map $\mu_k$ is
approximately unitarily equivalent to   $\id_B$.
But this follows from either \cite[Cor.~1.12]{TomsWinter:ssa} or Theorem~\ref{D-asu}.
\end{nnproof}
\en

\section{Proving the main result: The second approach}

This section is devoted to another proof of Theorem~\ref{C(X)-triviality}, which to some extent follows ideas of \cite{HirshbergRordamWinter:absorb-ssa}. Before turning to the actual proof, we outline our strategy.

\bn
\label{outline-wilhelms-approach}
We will obtain the isomorphism between $A$ and $\Ch(X) \otimes \Dh$ by constructing c.p.c.\ maps $\psi:A \to \Ch(X) \otimes \Dh$ and $\varphi: \Ch(X) \otimes \Dh \to A$ which implement an asymptotic  intertwining (in the sense of \cite{BlaKir:limits}) of  $A$ and $\Ch(X) \otimes \Dh$ as $\Ch(X)$-algebras (Proposition~\ref{C(X)-intertwining}).

Since $X$ is finite-dimensional, we may assume that $X \subset [0,1]^{n}$ for some $n \in \N$. We may then regard $A$ as a $\Ch(Y)$-algebra, where $Y = \pr_{1}(X) \subset [0,1]$ is the image of $X$ under the first coordinate projection. By an induction argument it will then be enough to prove that $A \cong \Ch(X) \otimes \Dh$ under the additional assumption that there is an isomorphism $\theta_{t}:A_{t} \cong \Ch(X_{t}) \otimes \Dh$ for each $t \in Y$, where $X_{t}= \pr^{-1}(\{t\})$ is the $(n-1)$-dimensional pre-image of $t$ in $X$. The  maps $\theta_{t}$ induce an embedding $\theta:A \to \prod_{t \in Y} \Ch(X_{t}) \otimes \Dh$.

The first problem then is that we do not know a priori how to choose these isomorphisms in a locally trivial manner. In other words, it is not clear whether we may assume that for any $t \in Y$ there are a neighborhood $Y^{t} \subset Y$ of $t$ and isomorphisms $\theta_{s}:A_{s} \cong \Ch(X_{s}) \otimes \Dh$, $s \in Y^{(t)}$, such that the induced map $\theta^{(t)}: A_{Y^{(t)}} \to \prod_{s \in Y^{(t)}} \Ch(X_{s}) \otimes \Dh$ actually maps $A_{Y^{(t)}}$ to $\Ch(X_{Y^{(t)}}) \otimes \Dh \subset \prod_{s \in Y^{(t)}} \Ch(X_{s}) \otimes \Dh$. However, this problem can be circumvented by the concept of `approximate local trivializations' $(\theta_{s}:A_{s} \stackrel{\cong}{\longrightarrow} \Ch(X_{s}))_{s \in Y^{(t)}}$ of A (introduced in Definition~\ref{d-local-approximate-trivialization}). The crucial property of such an approximate local trivialization is that it maps $A_{Y^{(t)}}$ to something \emph{close} to $\Ch(X_{Y^{(t)}}) \otimes \Dh \subset \prod_{s \in Y^{(t)}} \Ch(X_{s}) \otimes \Dh$. (This and the other technical properties of approximate local trivializations will be  controlled by the c.p.c.\ maps $\zeta^{(t)}: \Ch(X_{t}) \otimes \Dh \to \Ch(X_{Y^{(t)}}) \otimes \Dh$ and $\sigma^{(t)}: \Ch(X_{Y^{(t)}}) \otimes \Dh \to A_{Y^{(t)}}$ of \ref{d-local-approximate-trivialization}.) Using that $\Dh$ is strongly self-absorbing and $K_{1}$-injective, and involving an argument somewhat similar to Lemma~\ref{lemma:local-D-sections}, the existence of such local approximations will be established in Lemma~\ref{l-local-approximate-trivialization}.

Since $Y$ is compact and one-dimensional, we may then pick $0\le y_{0} < \ldots < y_{M}\le 1$ and cover $Y$ by closed sets $Y^{(y_{0})}, \ldots, Y^{(y_{M})}$ such that the intersections $Y^{(y_{i})} \cap Y^{(y_{i'})}$ are empty unless $|i-i'|\le1$. We set $I:= \{i \mid Y^{(y_{i})} \cap Y^{(y_{i+1})} \neq \emptyset\}$;  we may assume that for each $i \in I$ the set  $Y^{(y_{i})} \cap Y^{(y_{i'})}$ is a closed neighborhood of some $t_{i} \in Y$.

The next problem is that the local trivializations $\theta^{(y_{i})}$ over $Y^{(y_{i})}$ and $\theta^{(y_{i+1})}$ over $Y^{(y_{i+1})}$ do not necessarily agree over $Y^{(y_{i})} \cap Y^{(y_{i'})}$ (they need not even be close to each other in any way). However,  we may apply Lemma~\ref{l-local-approximate-trivialization} again to obtain a small interval $\widetilde{Y}^{(t_{i})} \subset Y^{(y_{i})} \cap Y^{(y_{i'})}$ around each $t_{i}$ (for each $i \in I$) and approximate local trivializations $\widetilde{\theta}^{(t_{i})}$ over $\widetilde{Y}^{(t_{i})}$. Again using that $\Dh$ is strongly self-absorbing and $K_{1}$-injective, this will allow us to modify the $\widetilde{\theta}^{(t_{i})}$ to obtain new approximate local trivializations (over  $\widetilde{Y}^{(t_{i})}$) which (approximately) intertwine $\theta^{(y_{i})}$ and $\theta^{(y_{i+1})}$ `along' $\widetilde{Y}^{(t_{i})}$.  Here, `along' means that the fibre maps $\theta^{(y_{i})}_{s}$ are continuously deformed into maps close to the fibre maps $\theta^{(y_{i+1})}_{s}$, where the parameter $s$ runs through the interval $\widetilde{Y}^{(t_{i})}$.

 By superposing the $\theta^{(y_{i})}$ and the (modified) $\widetilde{\theta}^{(t_{i})}$ we then obtain an approximate trivialization of $A$, which yields the desired approximate intertwining between $A$ and $\Ch(X) \otimes \Dh$.
\en

\bn
We start by  adapting  Elliott's approximate intertwining argument to  $\Ch(X)$-algebras as follows:

\label{C(X)-intertwining}
\begin{props}
Let $X$ be a compact metrizable space and let $A$ be a separable unital $\Ch(X)$-algebra with structure map $\mu$; let   $\Dh$ be a strongly self-absorbing $C^{*}$-algebra. Suppose that, for any $\varepsilon >0$ and for arbitrary compact subsets $\Fh \subset A$, $\Gh_{1} \subset \Ch(X)$ and $\Gh_{2} \subset \Dh$, there are c.p.c.\ maps
\[
\psi: A \to \Ch(X) \otimes \Dh
\]
and
\[
\varphi: \Ch(X) \otimes \Dh \to A
\]
satisfying
\begin{enumerate}
\item $\|\varphi \psi(a) - a\| < \varepsilon$ for $a \in \Fh$
\item $\|\varphi(f \otimes \be_{\Dh}) - \mu(f)\|< \varepsilon$ for $f \in \Gh_{1}$
\item $\|\psi \mu(f) - f \otimes \be_{\Dh}\| < \varepsilon$ for $f \in \Gh_{1}$
\item $\varphi$ is $\varepsilon$-multiplicative on $(\be_{X} \otimes \id_{\Dh})(\Gh_{2})$
\item $\psi$ is $\varepsilon$-multiplicative on $\Fh$.
\end{enumerate}
Then, $A$ and $\Ch(X) \otimes \Dh$ are isomorphic as $\Ch(X)$-algebras.
\end{props}

\begin{nproof}
From conditions (i)--(v) above in connection with Proposition \ref{prop;approx-D-au} it is straightforward to construct unitaries $u_{i} \in \Ch(X) \otimes \Dh$ and c.p.c.\ maps $\varphi_{i}: \Ch(X) \otimes \Dh \to A$ and $\psi_{i}: A \to \Ch(X) \otimes \Dh$, $i=1,2,\ldots$, such that the following diagram is an asymptotic intertwining in the sense of \cite[2.4]{BlaKir:limits}:
\begin{equation*}
\label{trivialization-diagram}
\xymatrix{
A  \ar[dr]^{\psi_{1}} \ar[rr]^{=} && A \ar[dr]^{\psi_{2}} \ar[rr]^{=} & &  \ldots \ar[dr] &  \\
& \Ch(X) \otimes \Dh \ar[ur]^{\varphi_{1}}  \ar[rr]^{\ad(\be_{X} \otimes u_{1})} && \Ch(X) \otimes \Dh  \ar[ur]^{\varphi_{2}} \ar[rr]^{\ad(\be_{X} \otimes u_{2})} && \ldots
}
\end{equation*}
Note that the inductive limit of the first row is just $A$, whereas the inductive limit of the second row is isomorphic to $\Ch(X) \otimes \Dh$. The $\varphi_{i}$ and $\psi_{i}$ then induce $*$-isomorphisms $\bar{\psi}:A \to \Ch(X) \otimes \Dh$ and $\bar{\varphi}:\Ch(X) \otimes \Dh \to A$ which are  mutual inverses (cf.\ the remark after Definition 2.4.1 of \cite{BlaKir:limits}). Using (ii), (iii) and the fact that the $u_{i}$ commute with $\Ch(X) \otimes \be_{\Dh}$, one can even assume that $\bar{\psi} \circ \mu = \id_{\Ch(X)} \otimes \be_{\Dh}$ and $\bar{\varphi} \circ  (\id_{\Ch(X)} \otimes \be_{\Dh}) = \mu$, which means that $\bar{\varphi}$ and $\bar{\psi}$ are isomorphisms of $\Ch(X)$-algebras.
\end{nproof}
\en

\bn
It will be convenient to note the following easy consequence of Proposition \ref{prop;approx-D-au} explicitly.

\label{au-bundle-maps}
\begin{lms}
Let $W$ be a compact metrizable space and $\Dh$ a strongly self-absorbing $K_{1}$-injective $C^{*}$-algebra. Then, for any compact subset $\be \in \Fh \subset \Ch(W) \otimes \Dh$ and $\gamma>0$ there are a compact subset $\Eh(\Fh,\gamma) \subset \Ch(W) \otimes \Dh$ and $0 < \delta(\Fh,\gamma) < \gamma/2$ such that the following holds:

If $K$ is a compact metrizable space and
\[
\sigma_{1},\sigma_{2} \colon \Ch(W) \otimes \Dh \to \Ch(K) \otimes \Dh
\]
are u.c.p.\ maps which map $\Ch(W)$ to $\Ch(K)$,  are $\delta(\Fh,\gamma)$-multiplicative on $\Eh(\Fh,\gamma)$ and agree on $\Ch(W)$, then there is a continuous path
\[
(u_{t})_{t \in [0,1]} \subset \Ch(K) \otimes \Dh
\]
of unitaries satisfying
\begin{equation}
\label{w1}
u_{0} = \be_{K} \otimes \be_{\Dh}
\mbox{ and } \|u_{1} \sigma_{1}(d) u_{1}^{*} - \sigma_{2}(d)\|< \gamma
\end{equation}
for $d \in \Fh \cdot \Fh$.

Moreover, we may assume that if $\Fh' $ is another compact subset containing $\Fh$ and $0 < \gamma' \le \gamma$, then  $  \Eh(\Fh,\gamma)\subset \Eh(\Fh',\gamma')$ and $\delta(\Fh',\gamma') < \delta(\Fh,\gamma)$.
\end{lms}

\begin{nproof}
It is straightforward to check that it suffices to prove the assertion when $\Fh$ is of the form
\[
\{f \otimes d \mid f \in \Fh_{0}, \, d \in \Fh_{1}\}
\]
where $\Fh_{0} \subset \Ch(W)$ and $\Fh_{1} \subset \Dh$ are compact subsets of normalized elements. We then apply Proposition \ref{prop;approx-D-au} (with $\be_{W} \otimes \Fh_{1} \cdot \Fh_{1}$ and $\sigma_{i} \circ (\be_{W} \otimes \id_{\Dh})$ in place of $\Fh$ and $\sigma_{i}$ for $i=0,1$) to obtain a compact subset $\Gh \subset \Dh$, $\delta>0$ and a continuous path $(u_{t})_{t\in [0,1]} \subset \Ch(K) \otimes \Dh$ of unitaries starting in $\be_{K} \otimes \be_{\Dh}$ and satisfying the assertion of \ref{prop;approx-D-au}.

Since, for $i=0,1$, $\sigma_{i}(\be_{W} \otimes \Dh)$ and $\sigma_{i}(\Ch(W) \otimes \be_{\Dh})$ commute, we see that $\sigma_{i}$ is the product map of the restrictions to $\Ch(W)$ and $\Dh$, respectively. In other words, we have $\sigma_{i}(f \otimes d) = \sigma_{i}(f \otimes \be_{\Dh}) \sigma_{i}(\be_{\Dh} \otimes d)$ for $f \in \Ch(W)$,  $d \in \Dh$ and $i=0,1$.

Setting $\Eh(\Fh,\gamma):= \be_{\Dh} \otimes \Gh$ and $\delta(\Fh,\gamma):=\delta$,  and using that each $u_{t}$ commutes with each $\sigma_{i}(f \otimes \be_{\Dh})$, we see that the assertion of \ref{prop;approx-D-au} yields \eqref{w1}.

We may clearly make $\delta(\Fh,\gamma)$ smaller and $\Eh(\Fh,\gamma)$ larger if necessary, from which the monotonicity statements follow.
\end{nproof}
\en

\bn
\label{A-Y-notation}
\begin{nots}
If $n \in \N$ and $X \subset [0,1]^{n}$ is a closed subset, we denote by $Y$ the image of $X$ under the first coordinate projection, $Y:= \pr_{1}(X) \subset [0,1]$. If $s \in Y$ (or if $V \subset Y$ is a closed subset), we set $X_{s}:= \pr_{1}^{-1}(\{s\})$ (or $X_{V}:= \pr_{1}^{-1}(V)$). If $A$ is a $\Ch(X)$-algebra, we write $A_{s}$ in place of $A_{X_{s}}$ (or $A_{V}$ in place of $A_{X_{V}}$); the fibre maps are then denoted by $\pi_{s}$ (or $\pi_{V}$). We may regard $A_{s}$  as a  $\Ch(X_{s})$-algebra and $A_{V}$ as a $\Ch(V)$-algebra in the obvious way.
\end{nots}
\en

\bn
\label{d-local-approximate-trivialization}
We will reduce our proof of Theorem \ref{C(X)-triviality} to the situation of \ref{A-Y-notation}; by induction we will be able to assume that each $A_{s}$, $s \in Y$, is isomorphic to $\Ch(X_{s}) \otimes \Dh$. However, it is not clear whether the isomorphisms can be chosen in a locally trivial manner. This technical problem is bypassed by the following  concept  of `approximate local trivializations':

\begin{dfs}
Let $n \in \N$, $X \subset [0,1]^{n}$ a closed subset and $A$ a unital $\Ch(X)$-algebra. Let $\Dh$ be another unital  $C^{*}$-algebra. Suppose each fibre of $A$ is isomorphic to $\Dh$ and that, for each $s \in Y$, we have $A_{s} \cong \Ch(X_{s}) \otimes \Dh$ as $\Ch(X)$-algebras.

Let $\eta>0$, $t \in Y$, an isomorphism $\theta \colon A_{t} \to \Ch(X_{t}) \otimes \Dh$ and  compact subsets $\be_{A} \in \Fh \subset A$, $\Gh \subset \Ch(X) \otimes \Dh$ and $\widehat{\Gh} \subset \Ch(X_{t}) \otimes \Dh$ be given.

By a $(\theta,\Fh,\Gh,\widehat{\Gh},\eta)$-trivialization of $A$ over the closed neighborhood $Y^{(t)}$ of $t$ we mean a  family of diagrams of $C^{*}$-algebras and u.c.p.\ maps
\begin{equation}
\label{approximate-trivialization-diagram}
\xymatrix{
&& A \ar[d]^{\pi_{Y^{(t)}}} & \\
&& A_{Y^{(t)}} \ar[r]^{\pi_{s}} \ar[d]^{\theta^{(t)}} & A_{s} \ar[d]^{\theta^{(t)}_{s}} \\
\Ch(X) \otimes \Dh \ar[r]^{\pi_{Y^{(t)}}} & \Ch(X_{Y^{(t)}}) \otimes \Dh \ar[r]^-{\iota^{(t)}} \ar[ur]^{\sigma^{(t)}} & \prod_{s' \in Y^{(t)}} \Ch(X_{s'}) \otimes \Dh \ar[r]^-{\pi_{s}} \ar[d]^{\pi_{t}}& \Ch(X_{s}) \otimes \Dh \\
&& \Ch(X_{t}) \otimes \Dh \ar[ul]^{\zeta^{(t)}}
}
\end{equation}
for each $s \in Y^{(t)}$  which satisfies the following properties:
\begin{enumerate}
\item all the $C^{*}$-algebras are $\Ch(X)$-algebras in the obvious way
\item  the maps $\pi_{s}$, $\pi_{Y^{(t)}}$ and $\iota^{(t)}$ are the obvious $\Ch(X)$-linear quotient and inclusion homomorphisms, respectively
\item  $\sigma^{(t)}$, $\theta^{(t)}$ and each $\theta^{(t)}_{s}$   are $\Ch(X)$-linear; $\zeta^{(t)}$ maps $\Ch(X_{t}) \otimes \be_{\Dh}$ to $\Ch(X_{Y^{(t)}}) \otimes \be_{\Dh}$
\item $\theta^{(t)}$ is a $*$-homomorphism and each $\theta^{(t)}_{s}$ is a $*$-isomorphism
\item $\theta^{(t)}_{t} = \theta$
\item the upper right rectangle commutes
\item any two paths starting and ending  in $\Ch(X_{t}) \otimes \Dh$ coincide (this in particular means that $\iota^{(t)} \circ \zeta^{(t)}$ lifts $\pi_{t}$ and that $\sigma^{(t)} \circ \zeta^{(t)}$ lifts $\theta^{(t)}_{t}$)
\item any two paths starting in $A$ and ending in the same algebra coincide up to $\eta$ on $\Fh \cdot \Fh$
\item any two paths starting in $\Ch(X) \otimes \Dh$ and ending in the same algebra coincide up to $\eta$ on $\Gh \cdot \Gh$
\item any two paths starting in $\Ch(X_{t}) \otimes \Dh$ and ending in the same algebra coincide up to $\eta$ on $\widehat{\Gh}$
\item any path starting in $\Ch(X_{t}) \otimes \Dh$ is $\delta^{(t)}$-bimultiplicative on $\Eh^{(t)}$ (cf.\ \ref{d-bimultiplicative}), where
\[
\Eh^{(t)}:= \Eh(\theta(\Fh) \cup \widehat{\Gh} \cup \pi_{t}(\Gh), \eta) \supset \Eh(\theta(\Fh),\eta) \cup \theta(\Fh) \cup \widehat{\Gh} \cup \pi_{t}(\Gh)
\]
and
\[
\delta^{(t)}:= \delta(\theta(\Fh) \cup \widehat{\Gh} \cup \pi_{t}(\Gh), \eta) < \eta/2
\]
come from Lemma \ref{au-bundle-maps}
\item any path starting in $\Ch(X) \otimes \Dh$ is $\eta$-multiplicative on $\Gh$.
\end{enumerate}
Here, by a path we mean a composition of maps from the diagram in which no map (except for the identity map) occurs more than once. We also regard the identity map of each algebra of the diagram as such a path.
\end{dfs}
\en

\bn
\label{restricting-local-approximate-trivialization}
\begin{rems}
It is obvious from the preceding definition that a $(\theta,\Fh,\Gh,\widehat{\Gh},\eta)$-trivialization of $A$ over  $Y^{(t)}$ may be restricted to smaller neighborhoods of $t$, i.e., if $\widetilde{Y}^{(t)} \subset Y^{(t)}$ is another closed neighborhood of $t$, then restricting the maps (and algebras) of diagram \eqref{approximate-trivialization-diagram} to $\widetilde{Y}^{(t)}$ yields a $(\theta,\Fh,\Gh,\widehat{\Gh},\eta)$-trivialization of $A$ over $\widetilde{Y}^{(t)}$.
\end{rems}
\en

\bn
\label{l-local-approximate-trivialization}
The following lemma establishes the existence of local approximate trivializations as in \ref{d-local-approximate-trivialization}.

\begin{lms}
Let $n \in \N$,  a closed subset $X \subset [0,1]^{n}$ and  a unital $\Ch(X)$-algebra $A$ be given. Let $\Dh$ be a $K_{1}$-injective strongly self-absorbing $C^{*}$-algebra. Suppose each fibre of $A$ is isomorphic to $\Dh$ and that, for each $s \in Y$, we have $A_{s} \cong \Ch(X_{s}) \otimes \Dh$ as $\Ch(X)$-algebras. Let $\eta>0$, $t \in Y$, an isomorphism $\theta \colon A_{t} \to \Ch(X_{t}) \otimes \Dh$ and  compact subsets $\be_{A} \in \Fh \subset A$, $\Gh \subset \Ch(X) \otimes \Dh$ and $\widehat{\Gh} \subset \Ch(X_{t}) \otimes \Dh$ be given. Then, there are a closed neighborhood $Y^{(t)} \subset Y$ of $t$ and a $(\theta,\Fh,\Gh,\widehat{\Gh},\eta)$-trivialization of $A$  over $Y^{(t)}$.
\end{lms}

\begin{nproof}
We first pick an arbitrary closed neighborhood $\widetilde{Y}^{(t)}$ of $t$. This yields maps $\pi_{\widetilde{Y}^{(t)}}$, $\widetilde{\iota}^{(t)}$ and $\widetilde{\pi}_{s}$ (for $s \in \widetilde{Y}^{(t)}$) as in diagram \eqref{approximate-trivialization-diagram}.  We then define maps $\widetilde{\zeta}^{(t)}$, $\widetilde{\sigma}^{(t)}$, $\widetilde{\theta}^{(t)}$ and $\widetilde{\theta}^{(t)}_{s}$ such that conditions (i)--(vi) of Definition \ref{d-local-approximate-trivialization} are satisfied. For a sufficiently small closed subset $Y^{(t)}$ of $\widetilde{Y}^{(t)}$ we  will then be able to modify the maps $\widetilde{\theta}^{(t)}$ and $\widetilde{\theta}^{(t)}_{s}$ to suitable new isomorphisms $\theta^{(t)}$ and $\theta^{(t)}_{s}$ for $s \in Y^{(t)}$. Restriction of our algebras and maps to  $Y^{(t)}$  will then yield the desired diagram \eqref{approximate-trivialization-diagram} with properties (i)--(xii).

So let us assume we have picked $\widetilde{Y}^{(t)}$. For $s \in \widetilde{Y}^{(t)}$ choose isomorphisms $\widetilde{\theta}^{(t)}_{s}: A_{s} \to \Ch(X_{s}) \otimes \Dh$; we take $\widetilde{\theta}^{(t)}_{t}$ to be $\theta$.

Next, we choose a c.p.c.\ lift $\widetilde{\zeta}^{(t)}: \Ch(X_{t}) \otimes \Dh \to \Ch(X_{\widetilde{Y}^{(t)}}) \otimes \Dh$ of $\pi_{t} \circ \widetilde{\iota}^{(t)}$; we may clearly assume that $\widetilde{\zeta}^{(t)}= \bar{\zeta}^{(t)} \otimes \id_{\Dh}$ for some u.p.c.\ map $\bar{\zeta}^{(t)}: \Ch(X_{t}) \to \Ch(X_{\widetilde{Y}^{(t)}})$.

Similarly, choose a u.p.c.\ lift $\bar{\sigma}^{(t)}: \Dh \to A_{\widetilde{Y}^{(t)}}$ of $\be_{X_{t}} \otimes \id_{\Dh}$ along $\widetilde{\pi}_{t}: A_{\widetilde{Y}^{(t)}} \to \Ch(X_{t}) \otimes \Dh$. Set $\widetilde{\sigma}^{(t)}:= \widetilde{\mu}^{(t)} \otimes \bar{\sigma}^{(t)}$, where $\widetilde{\mu}^{(t)}: \Ch(X_{\widetilde{Y}^{(t)}}) \to \Zh(A_{\widetilde{Y}^{(t)}})$ is the structure map of the $\Ch(\widetilde{Y}^{(t)})$-algebra $A_{\widetilde{Y}^{(t)}}$.

Define a compact subset of $\Ch(X_{\widetilde{Y}^{(t)}}) \otimes \Dh$ by
\[
\bar{\Fh}:= \widetilde{\zeta}^{(t)}(\widehat{\Gh}) \cup \widetilde{\pi}_{\widetilde{Y}^{(t)}}(\Gh \cdot \Gh) \cup  \widetilde{\zeta}^{(t)} \circ \theta   \circ \pi_{t}(\Fh \cdot \Fh).
\]
Apply Lemma \ref{au-bundle-maps} with $X_{\widetilde{Y}^{(t)}}$ in place of $W$ to obtain $\Eh(\bar{\Fh},\eta/2) \subset \Ch(X_{\widetilde{Y}^{(t)}}) \otimes \Dh$ and $0<\delta(\bar{\Fh},\eta/2)<\eta/4$.

Note that $\widetilde{\pi}^{(t)} \circ \widetilde{\sigma}^{(t)}$ is a $*$-homomorphism. It follows that there is a closed neighborhood $Y^{(t)} \subset \widetilde{Y}^{(t)}$ of $t$ such that $\widetilde{\pi}^{(s)} \circ \widetilde{\sigma}^{(t)}$ is $\delta(\bar{\Fh},\eta/2)$-multiplicative on $\Eh(\bar{\Fh},\eta/2)$ for all $s \in Y^{(t)}$. But then, by Lemma \ref{au-bundle-maps}, for each $s \in Y^{(t)}$ there is a unitary $\widetilde{u}^{(s)} \in \Ch(X_{s}) \otimes \Dh$ such that
\begin{equation}
\label{w2}
\| \widetilde{\pi}_{s} \circ \widetilde{\iota}^{(t)} (d)  - \widetilde{u}_{s}( \widetilde{\theta}^{(t)}_{s} \circ  \widetilde{\pi}^{(s)} \circ \widetilde{\sigma}^{(t)}(d) )\widetilde{u}_{s}^{*} \| < \eta/2
\end{equation}
for all $d \in \bar{\Fh}$.
Set $\theta^{(t)}_{t}:= \theta$; for $t \neq s \in Y^{(t)}$ we define
\[
\theta_{s}^{(t)} (\, . \,) := \widetilde{u}_{s}( \widetilde{\theta}^{(t)}_{s} \circ  \widetilde{\pi}^{(s)} (\, . \,) )\widetilde{u}_{s}^{*}.
\]
$\theta^{(t)}$ is defined accordingly.

Let $\pi_{Y^{(t)}}$, $\sigma^{(t)}$, $\zeta^{(t)}$, $\pi_{s}$ and $\iota^{(t)}$ denote the restrictions of $\widetilde{\pi}_{Y^{(t)}}$, $\widetilde{\sigma}^{(t)}$, $\widetilde{\zeta}^{(t)}$, $\widetilde{\pi}_{s}$ and $\widetilde{\iota}^{(t)}$, respectively, to $Y^{(t)}$.

We now have a diagram as \eqref{approximate-trivialization-diagram}. It is clear from our construction that (i)--(vii) of Definition \ref{d-local-approximate-trivialization} hold.

By making $Y^{(t)}$ smaller, if necessary, we may assume that (xi) and (xii) hold (using that $\pi_{t} \circ \sigma^{(t)}$ and $\pi_{t} \circ \zeta^{(t)}$ are exactly multiplicative and that $\sigma^{(t)}$ and $\zeta^{(t)}$ are lifts of these maps) and that
\begin{equation}
\label{w3}
\|\pi_{Y^{(t)}}(a) - \sigma^{(t)} \circ \zeta^{(t)} \circ \pi_{t} \circ \theta^{(t)}\circ \pi_{Y^{(t)}}(a)\|<\eta/2
\end{equation}
for $a \in \Fh \cdot \Fh$.

It follows from \eqref{w2} that $\pi_{s} \circ \iota^{(t)}$ and $\theta^{(t)}_{s} \circ \pi_{s} \circ \sigma^{(t)}$ coincide up to $\eta/2$ on $\bar{\Fh}$. Together with \eqref{w3} and commutativity of the upper right rectangle, this yields (viii), (ix) and (x) of Definition \ref{d-local-approximate-trivialization}.
\end{nproof}
\en

\bn
We are now prepared to give the second proof of our main result, following the outline in \ref{outline-wilhelms-approach}.

\begin{nnproof} (of Theorem~\ref{C(X)-triviality}):
Let $\varepsilon >0$ and  compact subsets $\Fh \subset A$, $\Gh_{1} \subset \Ch(X)$ and $\Gh_{2} \subset \Dh$ be given. We may assume that
\begin{equation}
\label{w23}
\be_{A} \in \Fh, \, \be_{X} \in \Gh_{1}, \, \be_{\Dh} \in \Gh_{2} \mbox{ and } \mu(\Gh_{1}) \subset \Fh.
\end{equation}
By Proposition~\ref{C(X)-intertwining}, it will be enough to  construct c.p.c.\ maps
\[
\psi: A \to \Ch(X) \otimes \Dh
\]
and
\[
\varphi: \Ch(X) \otimes \Dh \to A
\]
satisfying the conditions of \ref{C(X)-intertwining}.

By \cite[Theorem~V.3]{HurWal:Dim} we may assume that $X$ is a closed subset of $[0,1]^{n}$ for some $n \in \N$; we use the terminology of \ref{A-Y-notation}.

If $n=1$, then $X = Y$ and for each $t \in Y$ we have $X_{t} = \{t\}$, whence $A_{t} \cong \Dh \cong \Ch(X_{t}) \otimes \Dh$. Now if the theorem has been shown for closed subsets of $[0,1]^{n-1}$, it holds in particular with $A_{t}$ and  $X_{t}$ in place of $A$ and $X$ for any $t \in Y$. Therefore, for both the basis of induction and for the induction step we may assume that for each $t \in Y$ there is an isomorphism of $\Ch(X_{t})$-algebras  (hence also of $\Ch(X)$-algebras)
\[
\theta_{t} \colon A_{t} \to \Ch(X_{t}) \otimes \Dh.
\]
For convenience, set
\begin{equation}
\label{w18}
\eta:= \varepsilon/14
\end{equation}
and
\[
\Gh:= \{f \otimes d \mid f \in \Gh_{1}, d \in \Gh_{2}\}.
\]
Applying Lemma~\ref{l-local-approximate-trivialization} for each $t \in Y$ (with $\theta_{t}$ in place of $\theta$ and $\{0\} $ in place of $\widehat{\Gh}$) yields  compact neighborhoods $Y^{(t)} \subset Y$ of $t$ and $(\theta_{t},\Fh,\Gh,\widehat{\Gh},\eta)$-trivializations of $A$ over $Y^{(t)}$. Using that $Y$ is a compact subset of the one-dimensional space $[0,1]$ together with Remark~\ref{restricting-local-approximate-trivialization}, by making the $Y^{(t)}$ smaller, if necessary, it is straightforward to find
\begin{enumerate}
\item[a)] $M \in \N$, $0 = t_{0} \le y_{1} < t_{1} < y_{2} < \ldots < t_{M-1} < y_{M} \le t_{M} = 1$ and $\beta>0$
\item[b)] closed neighborhoods $Y^{(y_{i})} \subset Y$ of $y_{i}$, $i=1, \ldots, M$
\item[c)] $(\theta_{y_{i}},\Fh,\Gh,\widehat{\Gh},\eta)$-trivializations of $A$ over $Y^{(y_{i})}$ (as in diagram~\eqref{approximate-trivialization-diagram}, with $y_{i}$ in place of $t$), $i=1, \ldots, M$
\end{enumerate}
such that
\begin{enumerate}
\item[d)] $Y = \bigcup_{i=1}^{M} Y^{(y_{i})}$
\item[e)] each $Y^{(y_{i})}$ is the intersection of $Y$ with a closed interval
\item[f)] $Y^{(y_{i})} \cap Y^{(y_{i+2})} = \emptyset$, $i=1, \ldots, M-2$
\item[g)] $y_{i} \in [t_{i-1},t_{i}] \cap Y \subset Y^{(y_{i})}$  for $i=1, \ldots,M$
\item[h)] if $Y^{(y_{i})} \cap Y^{(y_{i+1})} \neq \emptyset$, we have $t_{i} \in [t_{i}-2 \beta,t_{i}+2 \beta] \cap Y \subset Y^{(y_{i})} \cap Y^{(y_{i+1})}$ 
\item[i)] if $Y^{(y_{i})} \cap Y^{(y_{i+1})} = \emptyset$, then $[t_{i}-2 \beta,t_{i}+2\beta] \cap Y = \emptyset$.
\end{enumerate}
We set
\[
I:= \{i \in \{1, \ldots,M-1\} \mid Y^{(y_{i})} \cap Y^{(y_{i+1})} \neq \emptyset \}
\]
and
\[
B^{(i)}:= Y^{(y_{i})} \cap Y^{(y_{i+1})}.
\]
For $i \in I$, we define u.c.p.\ maps
\[
\lambda^{(i)},\varrho^{(i)} \colon \Ch(X_{y_{i}}) \otimes \Dh \to \Ch(X_{t_{i}}) \otimes \Dh
\]
by
\begin{equation}
\label{w7}
\lambda^{(i)} :=        \theta^{(y_{i})}_{t_{i}} \pi_{t_{i}} \sigma^{(y_{i})} \zeta^{(y_{i})}  \mbox{ and } \varrho^{(i)} := \theta^{(y_{i+1})}_{t_{i}} \pi_{t_{i}} \sigma^{(y_{i})} \zeta^{(y_{i})}.
\end{equation}
The maps $\lambda^{(i)}$ and $\varrho^{(i)}$ send $\Ch(X_{y_{i}}) \otimes \be_{\Dh}$ to $\Ch(X_{t_{i}}) \otimes \be_{\Dh}$ and  coincide on $\Ch(X_{y_{i}}) \otimes \be_{\Dh}$ (as follows from \ref{d-local-approximate-trivialization}(iii)); they are $\delta^{(y_{i})}$-multiplicative on $\Eh(\theta^{(y_{i})}_{y_{i}}(\Fh),\eta) \subset \Eh^{(y_{i})}$ by \ref{d-local-approximate-trivialization}(xi) (the maps $\sigma^{(y_{i})} \circ \zeta^{(y_{i})}$ are -- and $\pi_{t_{i}}$, $\theta^{(y_{i})}_{t_{i}}$ and  $\theta^{(y_{i+1})}_{t_{i}}$ are exactly multiplicative). But then by Lemma~\ref{au-bundle-maps} there is a continuous path $(u^{(i)}_{t})_{t \in [0,1]}$ of unitaries in $\Ch(X_{t_{i}}) \otimes \Dh$ such that
\begin{equation}
\label{w8}
u_{0}^{(i)}= \be_{X_{t_{i}}} \otimes \be_{\Dh}
\end{equation}
and
\begin{equation}
\label{w9}
\|u_{1}^{(i)} \lambda^{(i)} \theta^{(y_{i})}_{y_{i}} \pi_{y_{i}}(a) (u^{(i)}_{1})^{*} - \varrho^{(i)} \theta^{(y_{i})}_{y_{i}} \pi_{y_{i}}(a)\|< \eta
\end{equation}
for $a \in \Fh \cdot \Fh$. We may as well regard the path $(u^{(i)}_{t})_{t \in [0,1]}$ as a unitary $u^{(i)} \in \Ch([0,1]) \otimes \Ch(X_{t_{i}}) \otimes \Dh$.

For each $i \in I$ we now apply Lemma~\ref{l-local-approximate-trivialization}   with $t_{i}$, $\theta^{(y_{i})}_{t_{i}} \pi_{y_{i}}$ and
\begin{equation}
\label{w5}
\widehat{\Gh}^{(i)}:= \{u^{(i)}_{t} \mid t \in [0,1]\}
\end{equation}
in place of $t$, $\theta$ and $\widehat{\Gh}$. We obtain for each $i \in I$ a closed neighborhood $\widetilde{Y}^{(t_{i})} \subset Y$ of $t_{i}$ and a $(\theta^{(y_{i})}_{t_{i}},\Fh,\Gh,\widehat{\Gh}^{(i)},\eta)$-trivialization of $A$ over $\widetilde{Y}^{(t_{i})}$. To make it easier to distinguish between these and the preceding approximate local trivializations, we denote the respective maps of diagram~\eqref{approximate-trivialization-diagram} by $\widetilde{\sigma}^{(t_{i})}$, $\widetilde{\zeta}^{(t_{i})}$, $\widetilde{\theta}^{(t_{i})}$ and $\widetilde{\theta}^{(t_{i})}_{s}$, respectively, and write $\widetilde{\Eh}^{(t_{i})}$ and $\widetilde{\delta}^{(t_{i})}$ in place of $\Eh(t)$ and $\delta^{(t)}$ (cf.~\ref{d-local-approximate-trivialization}(xi)).

Note that, in particular,
\begin{equation}
\label{w28}
\theta_{t_{i}}^{(y_{i})} = \widetilde{\theta}^{(t_{i})}_{t_{i}}
\end{equation}
by \ref{d-local-approximate-trivialization}(v).

Using  Remark~\ref{restricting-local-approximate-trivialization}, by making the $\widetilde{Y}^{(t_{i})}$ smaller if necessary, we may assume that  there is
\begin{equation}
\label{w15}
0<\beta'<\beta
\end{equation}
such that
\begin{enumerate}
\item[j)] $\widetilde{Y}^{(t_{i})} = [t_{i}-\beta',t_{i}+\beta'] \cap Y$ (this implies $\widetilde{Y}^{(t_{i})} \subset B^{(i)}$ by h))
\item[k)] for $t \in [t_{i}-2 \beta',t_{i}+2 \beta']$, $a \in \Fh$ and $d \in \Gh_{2}$ we have 
\begin{eqnarray*}
\lefteqn{\|\pi_{t}  ( \zeta^{(y_{i})} \theta_{y_{i}}^{(y_{i})} \pi_{y_{i}}(a) -  \widetilde{\zeta}^{(t_{i})} \widetilde{\theta}_{t_{i}}^{(t_{i})}\pi_{t_{i}}(a) )\| }\\
& \le &  \|\pi_{t_{i}}  ( \zeta^{(y_{i})} \theta_{y_{i}}^{(y_{i})} \pi_{y_{i}}(a) -  \widetilde{\zeta}^{(t_{i})} \widetilde{\theta}_{t_{i}}^{(t_{i})}\pi_{t_{i}}(a) )\| +  \eta
\end{eqnarray*}
and
\begin{eqnarray*}
\lefteqn{\|\pi_{t}  (\sigma^{(y_{i})} \zeta^{(y_{i})} \pi_{y_{i}}(\be_{X} \otimes d) - \widetilde{\sigma}^{(t_{i})} \pi_{\widetilde{Y}^{(t_{i})}}(\be_{X} \otimes d))\|}\\
& \le & \|\pi_{t_{i}}  (\sigma^{(y_{i})} \zeta^{(y_{i})} \pi_{y_{i}}(\be_{X} \otimes d) - \widetilde{\sigma}^{(t_{i})} \pi_{\widetilde{Y}^{(t_{i})}}(\be_{X} \otimes d))\| + \eta
\end{eqnarray*}
\item[l)] $\|\widetilde{\zeta}^{(t_{i})}(u_{1}^{(i)} \lambda^{(i)} \theta^{(y_{i})}_{y_{i}} \pi_{y_{i}}(a) (u^{(i)}_{1})^{*} - \varrho^{(i)} \theta^{(y_{i})}_{y_{i}} \pi_{y_{i}}(a) )\|< 2 \eta$ for $a \in \Fh \cdot \Fh$,
\end{enumerate}
where for k) we have used continuity of functions of the form ($\widetilde{Y}^{(t_{i})} \ni t \mapsto \|b(t)\|$) (cf.~\ref{upper-semicontinuity}); l) follows from \eqref{w9} and the fact that $\widetilde{\zeta}^{(t_{i})}$ lifts $\pi_{t_{i}}$ (by \ref{d-local-approximate-trivialization}(vii)).

For $i \in I$, define a function $h_{i} \in \Ch([t_{i}-2 \beta',t_{i} + 2\beta'])$ by
\begin{equation}
\label{w20}
h_{i}(t):= \left\{
\begin{array}{ll}
0 & \mbox{for } t \in [t_{i}-2 \beta',t_{i}-\beta']\\
1 & \mbox{for } t \in [t_{i}+ \beta',t_{i}+2\beta']\\
\frac{t-t_{i}+\beta'}{2\beta'} & \mbox{for } t \in [t_{i}- \beta',t_{i}+\beta']
\end{array}
\right.
\end{equation}
and an element
\[
\widetilde{u}^{(i)} \in \Ch(X_{\widetilde{Y}^{(t_{i})}}) \otimes \Dh
\]
by
\begin{equation}
\label{w4}
\pi_{t}(\widetilde{u}^{(i)}) := \pi_{t} \widetilde{\zeta}^{(t_{i})} (u_{h_{i}(t)}^{(i)})  \mbox{ for } t \in  \widetilde{Y}^{(t_{i})}.
\end{equation}
Define $(v_{t}^{(i)})_{t \in [0,1]} \in \Ch([0,1]) \otimes A_{\widetilde{Y}^{(t_{i})}}$ by
\begin{equation}
\label{w6}
v_{t}^{(i)}:= \widetilde{\sigma}^{(t_{i})} \widetilde{\zeta}^{(t_{i})} (u_{t}^{(i)})
\end{equation}
and
\[
\widetilde{v}^{(i)} := \Delta^{(i)}(v^{(i)}) \in A_{\widetilde{Y}^{(t_{i})}} ,
\]
where
\[
\Delta^{(i)} \colon \Ch([0,1]) \otimes A_{\widetilde{Y}^{(t_{i})}} \to A_{\widetilde{Y}^{(t_{i})}}
\]
is given by
\[
f \otimes a \mapsto (f \circ h) \cdot a;
\]
note that, with this definition,
\begin{equation}
\label{w14}
\pi_{t}(\widetilde{v}^{(i)}) = \pi_{t} \widetilde{\sigma}^{(t_{i})} \widetilde{\zeta}^{(t_{i})}(u^{(i)}_{h_{i}(t)})
\end{equation}
for $t \in \widetilde{Y}^{(t_{i})}$.

Next, choose positive functions $f_{1}, \ldots, f_{M},g_{1},\ldots,g_{M-1} \in \Ch(Y)$ such that
\begin{enumerate}
\item[m)] $\sum_{i=1}^{M}f_{i} + \sum_{i=1}^{M-1}g_{i} = \be_{Y}$
\item[n)] $\supp f_{i} \subset [t_{i-1}+\beta',t_{i}-\beta'] \cap Y \subset Y^{(y_{i})}  $ for $i=2, \ldots,M-1$,
\item[] $\supp f_{1} \subset [0,t_{1}-\beta'] \cap Y \subset Y^{(y_{1})}  $ and $\supp f_{M} \subset [t_{M-1}+\beta',1] \cap Y \subset Y^{(y_{M})}  $
\item[o)] $\supp g_{i} \subset [t_{i}-2\beta',t_{i} + 2\beta'] \cap Y \subset \widetilde{Y}^{(t_{i})}$ for $i \in I$ (we choose $g_{i}\equiv 0$ for $i \in \{1, \ldots,M-1\} \setminus I$).
\end{enumerate}
This is possible by a), h), i) and \eqref{w15}; note that
\begin{equation}
\label{w19}
f_{i}(y_{i'}) = g_{i}(t_{i'}) = \delta_{i,i'} \mbox{ and } g_{i}(y_{i'})=f_{i}(t_{i'}) = 0.
\end{equation}
We are finally prepared to define c.p.c.\ maps
\[
\psi \colon A \to \Ch(X) \otimes \Dh
\]
and
\[
\varphi \colon \Ch(X) \otimes \Dh \to A
\]
by
\begin{equation}
\label{w16}
\psi(\,.\,):= \sum_{i=1}^{M} f_{i} \cdot (\zeta^{(y_{i})}  \theta_{y_{i}}^{(y_{i})} \pi_{y_{i}})(\,.\,) + \sum_{i\in I} g_{i} \cdot (\ad_{\widetilde{u}^{(i)}} \widetilde{\zeta}^{(t_{i})}  \widetilde{\theta}^{(t_{i})}_{t_{i}}\pi_{t_{i}})(\, . \,)
\end{equation}
and
\begin{equation}
\label{w17}
\varphi(\,.\,):= \sum_{i=1}^{M} f_{i}\cdot (\sigma^{(y_{i})} \zeta^{(y_{i})} \pi_{y_{i}})(\,.\,) + \sum_{i\in I} g_{i}\cdot (\ad_{(\widetilde{v}^{(i)})^{*}}    \widetilde{\sigma}^{(t_{i})}  \pi_{\widetilde{Y}^{t_{i}}} )(\,.\,).
\end{equation}
It remains to check that $\psi$ and $\varphi$ satisfy  the conditions of Proposition~\ref{C(X)-intertwining}, i.e.,
\begin{enumerate}
\item[(i)] $\|\varphi \psi(a) - a\| < \varepsilon$ for $a \in \Fh$
\item[(ii)] $\|\varphi(f \otimes \be_{\Dh}) - \mu(f)\|< \varepsilon$ for $f \in \Gh_{1}$
\item[(iii)] $\|\psi \mu(f) - f \otimes \be_{\Dh}\| < \varepsilon$ for $f \in \Gh_{1}$
\item[(iv)] $\varphi$ is $\varepsilon$-multiplicative on $(\be_{X} \otimes \id_{\Dh})(\Gh_{2})$
\item[(v)] $\psi$ is $\varepsilon$-multiplicative on $\Fh$.
\end{enumerate}

For $a \in \Fh$, $i=1, \ldots, M-1$ and $t \in [t_{i-1}+2 \beta',t_{i}-2\beta'] \cap Y \subset Y^{(y_{i})}$ we have $f_{i}(t)=1$ and $g_{i}(t)=0$, whence
\begin{eqnarray*}
\|\pi_{t}(\varphi \psi(a)-a)\| & \stackrel{\eqref{w16},\eqref{w17}}{=} & \|\pi_{t} (\sigma^{(y_{i})} \zeta^{(y_{i})} \pi_{y_{i}}\zeta^{(y_{i})} \theta^{(y_{i})}_{y_{i}}  \pi_{y_{i}}(a) -a)\| \\
& \stackrel{\ref{d-local-approximate-trivialization}\mathrm{(vii)}}{\le} & \| \sigma^{(y_{i})} \zeta^{(y_{i})} \theta_{y_{i}}^{(y_{i})} \pi_{y_{i}}(a) -\pi_{Y^{(y_{i})}}(a)\| \\
& \stackrel{\ref{d-local-approximate-trivialization}\mathrm{(viii)}}{\le} & \eta\\
& \stackrel{\eqref{w18}}{<} &\varepsilon,
\end{eqnarray*}
and a similar estimate holds for $i = 1$ and $t \in [0,t_{1}-2\beta'] \cap Y$ and $i=M$ and $t \in [t_{M-1}+2\beta',1] \cap Y$.

For $a \in \Fh$, $i \in I$ and $t \in [t_{i} -\beta',t_{i}+\beta'] \cap Y \subset \widetilde{Y}^{(t_{i})}$ (cf.\ h) and j) above) we have $f_{i}(t)=0$ and $g_{i}(t)=1$ and therefore obtain
\begin{eqnarray*}
\lefteqn{\|\pi_{t} (\varphi \psi(a) - a)\|}\\
 & \stackrel{\eqref{w16},\eqref{w17}}{\le} & \|\pi_{t}( \ad_{(\widetilde{v}^{(i)})^{*}} \widetilde{\sigma}^{(t_{i})} \pi_{\widetilde{Y}^{(t_{i})}} \ad_{\widetilde{u}^{(i)}} \widetilde{\zeta}^{(t_{i})} \widetilde{\theta}_{t_{i}}^{(t_{i})} \pi_{t_{i}}(a) -a)\|\\
 & \stackrel{\eqref{w4}}{=} & \|\pi_{t}( \ad_{(\widetilde{v}^{(i)})^{*}} \widetilde{\sigma}^{(t_{i})}  \ad_{\widetilde{\zeta}^{(t_{i})}(u^{(i)}_{h_{i}(t)})} \widetilde{\zeta}^{(t_{i})} \widetilde{\theta}_{t_{i}}^{(t_{i})} \pi_{t_{i}}(a)) -a)\|\\
 & \stackrel{\ref{d-local-approximate-trivialization}\mathrm{(xi)},\eqref{w5}}{\le} & \|\pi_{t}( \ad_{(\widetilde{v}^{(i)})^{*}} \widetilde{\sigma}^{(t_{i})}  \widetilde{\zeta}^{(t_{i})} \ad_{u^{(i)}_{h_{i}(t)}}  \widetilde{\theta}_{t_{i}}^{(t_{i})} \pi_{t_{i}}(a)) -a)\| + \widetilde{\delta}^{(t_{i})}\\
 & \stackrel{\ref{d-local-approximate-trivialization}\mathrm{(xi)},\eqref{w5}}{\le} & \|\pi_{t}( \ad_{(\widetilde{v}^{(i)})^{*}} \ad_{\widetilde{\sigma}^{(t_{i})}\widetilde{\zeta}^{(t_{i})}(u^{(i)}_{h_{i}(t)})}  \widetilde{\sigma}^{(t_{i})}  \widetilde{\zeta}^{(t_{i})}  \widetilde{\theta}_{t_{i}}^{(t_{i})} \pi_{t_{i}}(a)) -a)\| + 2 \widetilde{\delta}^{(t_{i})}\\
  & \stackrel{\eqref{w6}}{=} & \| \pi_{t}( \ad_{(\widetilde{\sigma}^{(t_{i})} \widetilde{\zeta}^{(t_{i})} (u^{(i)}_{h_{i}(t)})^{*})} \ad_{\widetilde{\sigma}^{(t_{i})}\widetilde{\zeta}^{(t_{i})}(u^{(i)}_{h_{i}(t)})}\widetilde{\sigma}^{(t_{i})}  \widetilde{\zeta}^{(t_{i})} \widetilde{\theta}_{t_{i}}^{(t_{i})} \pi_{t_{i}}(a)) -a)\| \\
  && + 2 \widetilde{\delta}^{(t_{i})} \\
& = &  \| \pi_{t}( \ad_{(\widetilde{\sigma}^{(t_{i})} \widetilde{\zeta}^{(t_{i})} (u^{(i)}_{h_{i}(t)})^{*}) (\widetilde{\sigma}^{(t_{i})} \widetilde{\zeta}^{(t_{i})} (u^{(i)}_{h_{i}(t)}) )} \widetilde{\sigma}^{(t_{i})}  \widetilde{\zeta}^{(t_{i})} \widetilde{\theta}_{t_{i}}^{(t_{i})} \pi_{t_{i}}(a)) -a)\| \\
&& + 2 \widetilde{\delta}^{(t_{i})} \\
& \stackrel{\ref{d-local-approximate-trivialization}\mathrm{(xi)},\eqref{w5}}{<} & \| \pi_{t}(\widetilde{\sigma}^{(t_{i})} \widetilde{\zeta}^{(t_{i})} \widetilde{\theta}_{t_{i}}^{(t_{i})} \pi_{t_{i}}(a) - a)\|  + 4 \widetilde{\delta}^{(t_{i})} \\
& \stackrel{\ref{d-local-approximate-trivialization}\mathrm{(viii)}}{<} & \eta + 4 \widetilde{\delta}^{(t_{i})} \\
& \stackrel{\ref{d-local-approximate-trivialization}\mathrm{(xi)}}{<} & 3 \eta \\
& \stackrel{\eqref{w18}}{<} & \varepsilon.
\end{eqnarray*}
For $a \in \Fh$, $i=1, \ldots, M-1$ and $t \in [t_{i}-2 \beta', t_{i}-\beta'] \cap Y$ we have
\begin{eqnarray}
\label{w22}
\lefteqn{\|\pi_{t} (\varphi\psi(a)-a)\|} \nonumber \\
& = & \|\pi_{t} (f_{i} \cdot \sigma^{(y_{i})} \zeta^{(y_{i})} \theta^{(y_{i})}_{y_{i}} \pi_{y_{i}}(a)  \nonumber \\
& & + g_{i} \cdot f_{i} \cdot \ad_{(\widetilde{v}^{(i)})^{*}} \widetilde{\sigma}^{(t_{i})} \pi_{\widetilde{Y}^{(t_{i})}} \zeta^{(y_{i})} \theta_{y_{i}}^{(y_{i})} \pi_{y_{i}}(a)  \nonumber \\
&& + g_{i} \cdot g_{i} \cdot \ad_{(\widetilde{v}^{(i)})^{*}} \widetilde{\sigma}^{(t_{i})} \ad_{\widetilde{u}^{(i)}} \widetilde{\zeta}^{(t_{i})} \widetilde{\theta}_{t_{i}}^{(t_{i})} \pi_{t_{i}}(a) -a ) \|  \nonumber \\
& \stackrel{\eqref{w8},\eqref{w4},\eqref{w14},\eqref{w20}}{=} & \| \pi_{t} (f_{i} \cdot \sigma^{(y_{i})} \zeta^{(y_{i})} \theta^{(y_{i})}_{y_{i}} \pi_{y_{i}}(a)  + g_{i} \cdot f_{i} \cdot  \widetilde{\sigma}^{(t_{i})} \pi_{\widetilde{Y}^{(t_{i})}} \zeta^{(y_{i})} \theta_{y_{i}}^{(y_{i})} \pi_{y_{i}}(a)  \nonumber \\
&& + g_{i} \cdot g_{i} \cdot  \widetilde{\sigma}^{(t_{i})}  \widetilde{\zeta}^{(t_{i})} \widetilde{\theta}_{t_{i}}^{(t_{i})} \pi_{t_{i}}(a) -a ) \|  \nonumber \\
& \stackrel{\eqref{w21},\ref{d-local-approximate-trivialization}\mathrm{(viii)}}{\le} & \| f_{i}(t) \cdot \pi_{t}(a) + g_{i}(t) \cdot f_{i}(t) \cdot \pi_{t}(a) + g_{i}(t) \cdot g_{i}(t) \cdot  \pi_{t}(a) - \pi_{t}(a) \|  \nonumber \\
&&+ \eta + 3 \eta + \eta  \nonumber \\
& \stackrel{\mathrm{m)}}{=} & 5 \eta   \\
& \stackrel{\eqref{w18}}{<} & \varepsilon, \nonumber
\end{eqnarray}
where for the first equality we have used \eqref{w19},\eqref{w16},\eqref{w17},\ref{d-local-approximate-trivialization}(iii) and \ref{d-local-approximate-trivialization}(vii);  for the first inequality we have used
\begin{eqnarray}
\label{w21}
\lefteqn{\|\pi_{t} \widetilde{\sigma}^{(t_{i})} \pi_{\widetilde{Y}^{(t_{i})}} \zeta^{(y_{i})} \theta_{y_{i}}^{(y_{i})} \pi_{y_{i}}(a) - \pi_{t}(a) \|   } \nonumber \\
& \stackrel{\ref{d-local-approximate-trivialization}\mathrm{(viii)}}{\le} & \|\pi_{t} \widetilde{\sigma}^{(t_{i})} (\pi_{\widetilde{Y}^{(t_{i})}} \zeta^{(y_{i})} \theta_{y_{i}}^{(y_{i})} \pi_{y_{i}}(a) -  \widetilde{\zeta}^{(t_{i})} \widetilde{\theta}_{t_{i}}^{(t_{i})}\pi_{t_{i}}(a) )\| + \eta  \nonumber \\
& \le & \|\pi_{t}  ( \zeta^{(y_{i})} \theta_{y_{i}}^{(y_{i})} \pi_{y_{i}}(a) -  \widetilde{\zeta}^{(t_{i})} \widetilde{\theta}_{t_{i}}^{(t_{i})}\pi_{t_{i}}(a) )\| + \eta  \nonumber \\
& \stackrel{\mathrm{k)}}{\le} & \|\pi_{t_{i}}  ( \zeta^{(y_{i})} \theta_{y_{i}}^{(y_{i})} \pi_{y_{i}}(a) -  \widetilde{\zeta}^{(t_{i})} \widetilde{\theta}_{t_{i}}^{(t_{i})}\pi_{t_{i}}(a) )\| + 2 \eta  \nonumber \\
& \stackrel{\ref{d-local-approximate-trivialization}\mathrm{(vii),(viii)}}{<} & \|\theta_{t_{i}}^{(y_{i})} \pi_{t_{i}}(a) - \widetilde{\theta}_{t_{i}}^{(t_{i})} \pi_{t_{i}}(a)\| + 3 \eta  \nonumber \\
& = & 3 \eta  \nonumber \\
& < & \varepsilon.
\end{eqnarray}
A similar reasoning works for $t \in [t_{i}+\beta',t_{i}+2\beta'] \cap Y$. We have now checked condition (i) above.

For $f \in \Gh_{1}$, we have
\begin{eqnarray}
\label{w24}
\lefteqn{\|\psi \mu(f) - f \otimes \be_{\Dh}\|} \nonumber \\
& \stackrel{\eqref{w16}}{\le} & \|\sum_{i=1}^{M} f_{i} \cdot (\zeta^{(y_{i})} \theta^{(y_{i})}_{y_{i}} \pi_{y_{i}} \mu(f) - f \otimes \be_{\Dh}) \|  \nonumber \\
&  & + \|\sum_{i \in I} g_{i} \cdot (\ad_{\widetilde{u}^{(i)}} \widetilde{\zeta}_{t_{i}}^{(t_{i})} \pi_{t_{i}} \mu(f) - f \otimes \be_{\Dh}) \| \nonumber \\
& \stackrel{\eqref{w4}, \ref{d-local-approximate-trivialization}\mathrm{(xi)}}{\le} & \max_{i=1, \ldots,M} \|\pi_{Y^{(y_{i})}}(\zeta^{(y_{i})} \theta^{(y_{i})}_{y_{i}} \pi_{y_{i}}\mu(f) - f \otimes \be_{\Dh})\|  \nonumber \\
&& + \max_{i \in I} \| \ad_{\widetilde{u}^{(i)}} \widetilde{\zeta}^{(t_{i})} \widetilde{\theta}^{(t_{i})}_{t_{i}} \pi_{t_{i}} \mu(f) - \ad_{\widetilde{u}^{(i)}} (f \otimes \be_{\Dh})\| + \eta  \nonumber \\
& \stackrel{\ref{d-local-approximate-trivialization}\mathrm{(ii)}, \ref{d-local-approximate-trivialization}\mathrm{(ix)}}{<} & 3 \eta  \nonumber \\
& < & \varepsilon,
\end{eqnarray}
which yields (iii).

For (ii), note that
\begin{eqnarray*}
\|\varphi(f \otimes \be_{\Dh}) - \mu(f)\| & \le & \| \varphi(f \otimes \be_{\Dh}) - \varphi \psi\mu(f) \| + \|\varphi \psi \mu(f) - \mu(f)\| \\
& \stackrel{\eqref{w24},\eqref{w22},\eqref{w23}}{<} & 3 \eta + 6 \eta \\
& < & \varepsilon.
\end{eqnarray*}

Next, we check (iv). We have to show that
\begin{equation}
\|\pi_{t}(\varphi(\be_{X} \otimes d_{1} d_{2}) - \varphi(\be_{X} \otimes d_{1})\varphi(\be_{X} \otimes d_{2}) )\| < \varepsilon
\end{equation}
for any $d_{1},d_{2} \in \Gh_{2}$ and any $t \in Y$.

First, consider $t \in [t_{i-1}+2 \beta',t_{i} - 2 \beta'] \cap Y$ for $i=2, \ldots, M-1$. Then $f_{i}(t)=1$ and we have
\begin{eqnarray*}
\lefteqn{\|\pi_{t} (\varphi(\be_{X} \otimes d_{1} d_{2}) - \varphi(\be_{X} \otimes d_{1})\varphi(\be_{X} \otimes d_{2}))\|} \\
& = & \|\pi_{t} (\sigma^{(y_{i})} \zeta^{(y_{i})} \pi_{y_{i}} (\be_{X} \otimes d_{1} d_{2})\\
& & - \sigma^{(y_{i})} \zeta^{(y_{i})} \pi_{y_{i}}(\be_{X} \otimes d_{1})   \sigma^{(y_{i})} \zeta^{(y_{i})} \pi_{y_{i}} (\be_{X} \otimes d_{2}))\| \\
& \stackrel{\ref{d-local-approximate-trivialization}\mathrm{(xi)}}{<} & \delta^{(y_{i})} \\
& < & \varepsilon.
\end{eqnarray*}
A similar argument works for $t \in [0, t_{1}-2 \beta'] \cap Y$ and $t \in [t_{M-1}+ 2 \beta',1] \cap Y$.

If $t \in [t_{i}-\beta',t_{i}+\beta'] \cap Y$, then $g_{i}(t)=1$ and we have
\begin{eqnarray*}
\lefteqn{\|\pi_{t} (\varphi(\be_{X} \otimes d_{1}d_{2}) - \varphi(\be_{X} \otimes d_{1})\varphi(\be_{X} \otimes d_{2}))\|}\\
& = & \|\pi_{t}( \ad_{(\widetilde{v}^{(i)})^{*}} \widetilde{\sigma}^{(t_{i})} \pi_{\widetilde{Y}^{(t_{i})}}(\be_{X} \otimes d_{1}d_{2}) \\
&& - (\ad_{(\widetilde{v}^{(i)})^{*}} \widetilde{\sigma}^{(t_{i})} \pi_{\widetilde{Y}^{(t_{i})}}(\be_{X} \otimes d_{1}))(\ad_{(\widetilde{v}^{(i)})^{*}} \widetilde{\sigma}^{(t_{i})} \pi_{\widetilde{Y}^{(t_{i})}}(\be_{X} \otimes d_{2})))\| \\
& \stackrel{\eqref{w25},\eqref{w26}}{<} & \widetilde{\delta}^{(t_{i})} + \eta \\
& < & \varepsilon ,
\end{eqnarray*}
where we have used that
\begin{eqnarray}
\label{w25}
\lefteqn{ \|\pi_{t}(\widetilde{v}^{(i)}(\widetilde{v}^{(i)})^{*} - \be_{A_{\widetilde{Y}^{(t_{i})}}})\| } \nonumber \\
& = & \|\pi_{t}( \widetilde{\sigma}^{(t_{i})} \widetilde{\zeta}^{(t_{i})}(u_{h_{i}(t)}^{(i)}) \widetilde{\sigma}^{(t_{i})} \widetilde{\zeta}^{(t_{i})}(u_{h_{i}(t)}^{(i)})^{*} - \be_{\widetilde{Y}^{(t_{i})}} )\| \nonumber \\
& < & \widetilde{\delta}^{(t_{i})}
\end{eqnarray}
(by \ref{d-local-approximate-trivialization}(xi) and since the $u^{(i)}_{h_{i}(t)}$ are unitaries) and that
\begin{equation}
\label{w26}
 \| \widetilde{\sigma}^{(t_{i})}\pi_{\widetilde{Y}^{(t_{i})}}(\be_{X} \otimes d_{1} d_{2}) - \widetilde{\sigma}^{(t_{i})}\pi_{\widetilde{Y}^{(t_{i})}}(\be_{X} \otimes d_{1}) \widetilde{\sigma}^{(t_{i})}\pi_{\widetilde{Y}^{(t_{i})}}(\be_{X} \otimes d_{2}) \|
 \stackrel{\ref{d-local-approximate-trivialization}\mathrm{(xii)}}{<}  \eta.
\end{equation}

Now if $t \in [t_{i}-2\beta',t_{i}-\beta'] \cap Y$, we have
\begin{equation}
\label{w27}
\pi_{t}(\widetilde{v}^{(i)}) \stackrel{\eqref{w20},\eqref{w14}}{=} \pi_{t}(\be_{A})
\end{equation}
and obtain
\begin{eqnarray*}
\lefteqn{\|\pi_{t} (\varphi(\be_{X} \otimes d_{1}d_{2}) - \varphi(\be_{X} \otimes d_{1})\varphi(\be_{X} \otimes d_{2}))\|}\\
& \stackrel{\eqref{w16},\eqref{w17},\eqref{w27}}{=} & \| f_{i}(t)  (f_{i}(t) + g_{i}(t)) \cdot \pi_{t} \sigma^{(y_{i})} \zeta^{(y_{i})} \pi_{y_{i}}(\be_{X} \otimes d_{1}d_{2}) \\
&& + g_{i}(t) (f_{i}(t) + g_{i}(t)) \cdot \pi_{t} \widetilde{\sigma}^{(t_{i})} \pi_{\widetilde{Y}^{(t_{i})}} (\be_{X} \otimes d_{1} d_{2}) \\
&& - f_{i}(t)^{2} \cdot \pi_{t}( \sigma^{(y_{i})} \zeta^{(y_{i})} \pi_{y_{i}}(\be_{X} \otimes d_{1})  \sigma^{(y_{i})} \zeta^{(y_{i})} \pi_{y_{i}}(\be_{X} \otimes d_{2})) \\
&& - g_{i}(t)^{2} \cdot \pi_{t} (\widetilde{\sigma}^{(t_{i})} \pi_{\widetilde{Y}^{(t_{i})}} (\be_{X} \otimes d_{1})  \widetilde{\sigma}^{(t_{i})} \pi_{\widetilde{Y}^{(t_{i})}} (\be_{X} \otimes  d_{2})) \\
&& - f_{i}(t)  g_{i}(t) \cdot \pi_{t} ( \sigma^{(y_{i})} \zeta^{(y_{i})} \pi_{y_{i}}(\be_{X} \otimes d_{1}) \widetilde{\sigma}^{(t_{i})} \pi_{\widetilde{Y}^{(t_{i})}} (\be_{X} \otimes  d_{2}) )\\
&& - f_{i}(t)  g_{i}(t) \cdot \pi_{t} (\widetilde{\sigma}^{(t_{i})} \pi_{\widetilde{Y}^{(t_{i})}} (\be_{X} \otimes d_{1} ) \sigma^{(y_{i})} \zeta^{(y_{i})} \pi_{y_{i}}(\be_{X} \otimes d_{2})) \| \\
& \stackrel{\ref{d-local-approximate-trivialization}\mathrm{(xii)}}{\le} & f_{i}(t)^{2} \cdot \|\pi_{t} (\sigma^{(y_{i})} \zeta^{(y_{i})} \pi_{y_{i}} (\be_{X} \otimes d_{1} d_{2}) \\
&& - \sigma^{(y_{i})} \zeta^{(y_{i})} \pi_{y_{i}} (\be_{X} \otimes d_{1} )\sigma^{(y_{i})} \zeta^{(y_{i})} \pi_{y_{i}} (\be_{X} \otimes d_{2}))\| \\
&& + g_{i}(t)^{2} \cdot \|\pi_{t} (\widetilde{\sigma}^{(t_{i})}  \pi_{\widetilde{Y}^{t_{i}}} (\be_{X} \otimes d_{1} d_{2}) \\
&& - \widetilde{\sigma}^{(t_{i})}  \pi_{\widetilde{Y}^{(t_{i})}} (\be_{X} \otimes d_{1} ) \widetilde{\sigma}^{(t_{i})}  \pi_{\widetilde{Y}^{(t_{i})}} (\be_{X} \otimes d_{2}))\| \\
& & + f_{i}(t) g_{i}(t) \cdot (\| \pi_{t}(\sigma^{(y_{i})} \zeta^{(y_{i})} \pi_{y_{i}}(\be_{X} \otimes d_{1}) \sigma^{(y_{i})} \zeta^{(y_{i})} \pi_{y_{i}}(\be_{X} \otimes d_{2})  \\
&& +   \widetilde{\sigma}^{(t_{i})} \pi_{\widetilde{Y}^{(t_{i})}} (\be_{X} \otimes d_{1}) \widetilde{\sigma}^{(t_{i})} \pi_{\widetilde{Y}^{(t_{i})}} (\be_{X} \otimes d_{2})  \\
&& -    \sigma^{(y_{i})} \zeta^{(y_{i})} \pi_{y_{i}}(\be_{X} \otimes d_{1})  \widetilde{\sigma}^{(t_{i})} \pi_{\widetilde{Y}^{(t_{i})}} (\be_{X} \otimes d_{2})\\
&& -   \widetilde{\sigma}^{(t_{i})} \pi_{\widetilde{Y}^{(t_{i})}} (\be_{X} \otimes d_{1})  \sigma^{(y_{i})} \zeta^{(y_{i})} \pi_{y_{i}}(\be_{X} \otimes d_{2})   ) \| + 2 \eta   )\\
& \stackrel{\ref{d-local-approximate-trivialization}\mathrm{(xi)},\ref{d-local-approximate-trivialization}\mathrm{(xii)}}{<} & (f_{i}(t)^{2} + g_{i}(t)^{2}) \cdot \eta \\
&& + 2 f_{i}(t) g_{i}(t) \cdot \|\pi_{t}  (\sigma^{(y_{i})} \zeta^{(y_{i})} \pi_{y_{i}}(\be_{X} \otimes d_{1}) - \widetilde{\sigma}^{(t_{i})} \pi_{\widetilde{Y}^{(t_{i})}}(\be_{X} \otimes d_{1}))\| \\
& \stackrel{\mathrm{k)}}{\le}&  (f_{i}(t)^{2} + g_{i}(t)^{2}) \cdot \eta  + 2 f_{i}(t) g_{i}(t) \cdot \eta\\
&&  + \|\pi_{t_{i}}  (\sigma^{(y_{i})} \zeta^{(y_{i})} \pi_{y_{i}}(\be_{X} \otimes d_{1}) - \widetilde{\sigma}^{(t_{i})} \pi_{\widetilde{Y}^{(t_{i})}}(\be_{X} \otimes d_{1}))\| \\
& \stackrel{\ref{d-local-approximate-trivialization}\mathrm{(v)}}{=} &  \|\theta^{(y_{i})}_{t_{i}} \pi_{t_{i}} \sigma^{(y_{i})} \zeta^{(y_{i})} \pi_{y_{i}} (\be_{X} \otimes d_{1}) \\
&& - \widetilde{\theta}^{(t_{i})}_{t_{i}} \pi_{t_{i}} \widetilde{\sigma}^{(t_{i})} \pi_{\widetilde{Y}^{(t_{i})}}(\be_{X} \otimes d_{1}) \| + \eta \\
& \stackrel{\ref{d-local-approximate-trivialization}\mathrm{(ix)}}{\le} &  \| \pi_{t_{i}} (\be_{X} \otimes d_{1}) - \pi_{t_{i}}(\be_{X} \otimes d_{1})\| + 3 \eta \\
& = & 3 \eta \\
& < & \varepsilon.
\end{eqnarray*}
A similar reasoning works for $t \in [t_{i}+\beta',t_{i}+2 \beta'] \cap Y$, whence (iv) above  holds.

Property (v) is checked in a similar manner: we have to show that
\begin{equation}
\label{w10}
\|\pi_{t}(\psi(a_{1}a_{2}) - \psi(a_{1}) \psi(a_{2}))\| < \varepsilon
\end{equation}
for any $a_{1},a_{2} \in \Fh$ and any $t \in Y$. If $f_{i}(t)=1$ or $g_{i}(t)=1$, this follows immediately from \ref{d-local-approximate-trivialization}(xi) and our construction. Thus, it remains to check \eqref{w10} for $t \in [t_{i}-2 \beta',t_{i} - \beta'] \cap Y$, $i \in  I$ (the case $t \in [t_{i}+\beta',t_{i}+2\beta']$ runs parallel).

Let us first observe that for $a \in \Fh \cdot \Fh$
\begin{eqnarray}
\label{w13}
\lefteqn{\| \pi_{t}(\widetilde{\zeta}^{(t_{i})}  \widetilde{\theta}^{(t_{i})}_{t_{i}} \pi_{t_{i}}(a) - \zeta^{(y_{i})} \theta^{(y_{i})}_{y_{i}} \pi_{y_{i}}(a))\|} \nonumber \\
& \stackrel{\mathrm{k)}}{\le} & \| \pi_{t_{i}}(\widetilde{\zeta}^{(t_{i})} \widetilde{\theta}^{(t_{i})}_{t_{i}} \pi_{t_{i}}(a) - \zeta^{(y_{i})} \theta^{(y_{i})}_{y_{i}} \pi_{y_{i}}(a))\| +  \eta \nonumber \\
& \stackrel{\ref{d-local-approximate-trivialization}\mathrm{(vii)},\ref{d-local-approximate-trivialization}\mathrm{(viii)}}{<} & \|\widetilde{\theta}^{(t_{i})}_{t_{i}} \pi_{t_{i}}(a) - \theta^{(y_{i})}_{t_{i}} \pi_{y_{i}}(a)\| + 2 \eta \nonumber \\
& \stackrel{\eqref{w28}}{=} & 2 \eta.
\end{eqnarray}
We are now ready to compute
\begin{eqnarray*}
\lefteqn{\|\pi_{t}(\psi(a_{1}a_{2}) - \psi(a_{1}) \psi(a_{2}))\|} \\
& \stackrel{\eqref{w16},\eqref{w8},\eqref{w20},\eqref{w4}}{=} & \| f_{i}(t) (f_{i}(t)+g_{i}(t)) \cdot \pi_{t} \zeta^{(y_{i})} \theta^{(y_{i})}_{y_{i}} \pi_{y_{i}}(a_{1}a_{2}) \\
& & + g_{i}(t)(f_{i}(t) +g_{i}(t)) \cdot \pi_{t}  \widetilde{\zeta}^{(t_{i})} \widetilde{\theta}^{(t_{i})}_{t_{i}} \pi_{t_{i}}(a_{1}a_{2}) \\
& & - f_{i}(t)^{2} \cdot \pi_{t}( \zeta^{(y_{i})} \theta^{(y_{i})}_{y_{i}} \pi_{y_{i}}(a_{1}) \cdot \zeta^{(y_{i})} \theta^{(y_{i})}_{y_{i}} \pi_{y_{i}}(a_{2})) \\
&& - g_{i}(t)^{2} \cdot \pi_{t} (\widetilde{\zeta}^{(t_{i})} \widetilde{\theta}^{(t_{i})}_{t_{i}} \pi_{t_{i}}(a_{1}) \cdot \widetilde{\zeta}^{(t_{i})} \widetilde{\theta}^{(t_{i})}_{t_{i}} \pi_{t_{i}}(a_{2})) \\
&& - f_{i}(t) g_{i}(t) \cdot \pi_{t}(\zeta^{(y_{i})} \theta^{(y_{i})}_{y_{i}} \pi_{y_{i}}(a_{1})  \cdot \widetilde{\zeta}^{(t_{i})} \widetilde{\theta}^{(t_{i})}_{t_{i}} \pi_{t_{i}}(a_{2})  ) \\
&& - f_{i}(t) g_{i}(t) \cdot \pi_{t}(  \widetilde{\zeta}^{(t_{i})} \widetilde{\theta}^{(t_{i})}_{t_{i}} \pi_{t_{i}}(a_{1}) \cdot \zeta^{(y_{i})} \theta^{(y_{i})}_{y_{i}} \pi_{y_{i}}(a_{2}) ) \| \\
& \stackrel{\ref{d-local-approximate-trivialization}\mathrm{(xi)}}{\le} & \| f_{i}(t) g_{i}(t) \cdot \pi_{t} \zeta^{(y_{i})} \theta^{(y_{i})}_{y_{i}} \pi_{y_{i}}(a_{1}a_{2}) \\
& & + f_{i}(t)g_{i}(t)  \cdot \pi_{t}  \widetilde{\zeta}^{(t_{i})} \widetilde{\theta}^{(t_{i})}_{t_{i}} \pi_{t_{i}}(a_{1}a_{2}) \\
& & - f_{i}(t) g_{i}(t) \cdot \pi_{t}(\zeta^{(y_{i})} \theta^{(y_{i})}_{y_{i}} \pi_{y_{i}}(a_{1})  \cdot \widetilde{\zeta}^{(t_{i})} \widetilde{\theta}^{(t_{i})}_{t_{i}} \pi_{t_{i}}(a_{2})  ) \\
& & - f_{i}(t) g_{i}(t) \cdot \pi_{t}(  \widetilde{\zeta}^{(t_{i})} \widetilde{\theta}^{(t_{i})}_{t_{i}} \pi_{t_{i}}(a_{1}) \cdot \zeta^{(y_{i})} \theta^{(y_{i})}_{y_{i}} \pi_{y_{i}}(a_{2}) ) \| \\
& & + \delta^{(y_{i})} + \widetilde{\delta}^{(t_{i})} \\
& \stackrel{\eqref{w13}}{<} & 2 f_{i}(t) g_{i}(t) \cdot \| \pi_{t}( \zeta^{(y_{i})} \theta^{(y_{i})}_{y_{i}} \pi_{y_{i}}(a_{1}a_{2}) \\
& & - \zeta^{(y_{i})} \theta^{(y_{i})}_{y_{i}} \pi_{y_{i}}(a_{1}) \cdot  \zeta^{(y_{i})} \theta^{(y_{i})}_{y_{i}} \pi_{y_{i}}(a_{2})\| \\
& & + \delta^{(y_{i})} +  \widetilde{\delta}^{(t_{i})} + 4 \eta\\
& \stackrel{\ref{d-local-approximate-trivialization}\mathrm{(xi)}}{<} & 2 \delta^{(y_{i})} + \widetilde{\delta}^{(t_{i})} + 4 \eta\\
& < & \varepsilon
\end{eqnarray*}
for $a_{1},a_{2} \in \Fh$ and $t \in [t_{i}-2\beta',t_{i}-\beta'] \cap Y$. This completes the proof.
\end{nnproof}
\en


\bibliographystyle{amsplain}

\end{document}